\documentclass[hidelinks,onefignum,onetabnum]{siamart251216}

\usepackage{amsfonts,amssymb,mathtools}
\usepackage{graphicx}
\usepackage{booktabs}
\usepackage{algorithm}
\usepackage{algpseudocode}
\usepackage{placeins}
\usepackage[shortlabels]{enumitem}
\setlist[itemize]{leftmargin=*}

\ifpdf
  \DeclareGraphicsExtensions{.eps,.pdf,.png,.jpg}
\else
  \DeclareGraphicsExtensions{.eps}
\fi

\newsiamremark{remark}{Remark}
\crefname{remark}{Remark}{Remarks}
\Crefname{remark}{Remark}{Remarks}

\headers{Galerkin time stepping with linearized neural networks}{Youngkyu Lee}

\title{Galerkin Time Stepping with Linearized Neural Networks for Stiff Linear Evolution
Equations\thanks{Submitted to the editors \today.
\funding{This work was supported by a grant from Kyung Hee University in 2026 (KHU-20263365).}}}

\author{Youngkyu Lee\thanks{Department of Mathematics, Kyung Hee University, Seoul 02447, Republic of Korea
  (\email{youngkyu\_lee@khu.ac.kr}).}}

\ifpdf
\hypersetup{pdftitle={Galerkin Time Stepping with Linearized Neural Networks for Stiff Linear Evolution Equations},
  pdfauthor={Youngkyu Lee}}
\fi

\usepackage{amsopn}
\newcommand{\R}{\mathbb{R}}
\newcommand{\C}{\mathbb{C}}
\newcommand{\Pq}{\mathbb{P}}
\newcommand{\hI}{\hat I}
\newcommand{\hS}{\hat S}
\newcommand{\ip}[2]{(#1,#2)}
\newcommand{\norm}[1]{\lVert #1\rVert}
\newcommand{\anorm}[1]{\lvert #1\rvert_a}
\newcommand{\dd}{\,\mathrm{d}}
\DeclareMathOperator{\spanop}{span}
\DeclareMathOperator{\diag}{diag}
\renewcommand{\Re}{\operatorname{Re}}

\newcommand{\tablebodyfont}{\small}
\newsavebox{\tablebox}
\newcommand{\usetable}{\typeout{TABLEWIDTH \the\wd\tablebox\space of \the\textwidth}%
  \begin{tabular}{@{}c@{}}\ifdim\wd\tablebox>\textwidth\resizebox{\textwidth}{!}{\usebox{\tablebox}}%
  \else\usebox{\tablebox}\fi\end{tabular}}

\newcommand{\univr}{118}
\newcommand{\certELeg}{\ensuremath{1.1{\times}10^{-2}}}
\newcommand{\certESphere}{\ensuremath{2.2{\times}10^{-4}}}
\newcommand{\certERandom}{\ensuremath{3.4{\times}10^{-3}}}
\newcommand{\dgMinReMin}{0.65}
\newcommand{\univEone}{\ensuremath{5.6{\times}10^{-11}}}
\newcommand{\tailOthersDev}{\ensuremath{1.3{\times}10^{-9}}}
\newcommand{\stabGalMax}{\ensuremath{1.3{\times}10^{-15}}}
\newcommand{\stabGradedMax}{\ensuremath{5.2{\times}10^{-7}}}
\newcommand{\gradedNormA}{\ensuremath{6.7{\times}10^{11}}}
\newcommand{\dgRinfMax}{\ensuremath{6{\times}10^{-11}}}
\newcommand{\subUnknowns}{119}
\newcommand{\univInside}{\ensuremath{9.2{\times}10^{-10}}}
\newcommand{\ablRandomInside}{\ensuremath{4.1{\times}10^{-8}}}
\newcommand{\hpOptInside}{\ensuremath{2.0{\times}10^{-8}}}
\newcommand{\hpOptMid}{\ensuremath{1.1{\times}10^{-13}}}
\newcommand{\hpOptEnd}{\ensuremath{<10^{-13}}}
\newcommand{\hpOptRatio}{0.25}
\newcommand{\hpOptLevels}{9}
\newcommand{\hpOptQa}{3}
\newcommand{\hpOptQb}{21}
\newcommand{\hpSearched}{336}
\newcommand{\heatN}{16{,}129}
\newcommand{\heatRefX}{\ensuremath{1.7{\times}10^{-11}}}
\newcommand{\heatMaxbestRadau}{\ensuremath{5.1{\times}10^{-11}}}
\newcommand{\heatMaxbesttimeRadau}{15.49}
\newcommand{\heatMaxDGeight}{\ensuremath{3.8{\times}10^{-12}}}
\newcommand{\heatTimeDGeight}{12.23}
\newcommand{\heatSetupDG}{6.8}
\newcommand{\heatMaxAdeight}{\ensuremath{3.8{\times}10^{-12}}}
\newcommand{\heatTimeAdeight}{11.26}
\newcommand{\heatFactAd}{159}
\newcommand{\heatMaxAdPeight}{\ensuremath{3.8{\times}10^{-12}}}
\newcommand{\heatTimeAdPeight}{8.87}
\newcommand{\heatFactAdP}{134}
\newcommand{\heatMaxHPctwelve}{\ensuremath{2.1{\times}10^{-12}}}
\newcommand{\heatTimeHPctwelve}{6.91}
\newcommand{\heatFactHPctwelve}{105}
\newcommand{\heatMaxHPcfourteen}{\ensuremath{2.1{\times}10^{-13}}}
\newcommand{\heatTimeHPcfourteen}{14.90}
\newcommand{\heatFactHPcfourteen}{225}
\newcommand{\heatESDIRKeight}{\ensuremath{3.1{\times}10^{-7}}}
\newcommand{\heatESDIRKeighttime}{1.35}
\newcommand{\heatSpaceErrFirst}{\ensuremath{6.0{\times}10^{-3}}}
\newcommand{\heatSpaceErr}{\ensuremath{1.8{\times}10^{-4}}}
\newcommand{\meshNmax}{65{,}025}
\newcommand{\meshNmin}{961}
\newcommand{\meshGradedErr}{\ensuremath{1.0{\times}10^{-11}}}
\newcommand{\meshRefExact}{\ensuremath{3.1{\times}10^{-13}}}
\newcommand{\meshLamFirst}{\ensuremath{1.3{\times}10^{3}}}
\newcommand{\meshLamLast}{\ensuremath{8.5{\times}10^{4}}}
\newcommand{\refQ}{28}
\newcommand{\refL}{16}
\newcommand{\aprQ}{20}
\newcommand{\aprL}{12}
\newcommand{\timeGradedMax}{\ensuremath{8.9{\times}10^{-12}}}
\newcommand{\timeGradedUnk}{952}
\newcommand{\timeAdMax}{\ensuremath{8.9{\times}10^{-12}}}
\newcommand{\timeAdUnk}{399}
\newcommand{\timeAdPMax}{\ensuremath{8.9{\times}10^{-12}}}
\newcommand{\timeAdPUnk}{224}
\newcommand{\timeHPcMax}{\ensuremath{2.5{\times}10^{-13}}}
\newcommand{\timeHPcUnk}{330}
\newcommand{\timeAprMax}{\ensuremath{4.6{\times}10^{-11}}}
\newcommand{\timeAprUnk}{303}
\newcommand{\timeELMGrMax}{\ensuremath{4.5{\times}10^{-10}}}
\newcommand{\timeLegFirst}{\ensuremath{2.3{\times}10^{-2}}}
\newcommand{\timeRadauFirst}{\ensuremath{2.4{\times}10^{-2}}}
\newcommand{\meshLegFirst}{\ensuremath{3.1{\times}10^{-13}}}
\newcommand{\meshLegLast}{\ensuremath{5.3{\times}10^{-4}}}
\newcommand{\meshSphereLast}{\ensuremath{6.6{\times}10^{-3}}}
\newcommand{\meshRandomLast}{\ensuremath{7.6{\times}10^{-3}}}
\newcommand{\meshAprFirst}{\ensuremath{2.6{\times}10^{-12}}}
\newcommand{\meshAprLast}{\ensuremath{2.9{\times}10^{-10}}}
\newcommand{\meshAdPErr}{\ensuremath{1.0{\times}10^{-11}}}
\newcommand{\meshHPcErr}{\ensuremath{3.1{\times}10^{-13}}}
\newcommand{\refCondS}{\ensuremath{6.9{\times}10^{11}}}
\newcommand{\refErrSchurZero}{\ensuremath{3.1{\times}10^{-10}}}
\newcommand{\refErrSchurOne}{\ensuremath{7.1{\times}10^{-12}}}
\newcommand{\refErrSchurTwo}{\ensuremath{7.1{\times}10^{-12}}}
\newcommand{\refErrDiagFour}{\ensuremath{3.2{\times}10^{-7}}}
\newcommand{\refErrDiagSix}{\ensuremath{4.5{\times}10^{-11}}}
\newcommand{\refErrDiagEight}{\ensuremath{5.0{\times}10^{-11}}}
\newcommand{\refErrAprZero}{\ensuremath{3.4{\times}10^{-11}}}
\newcommand{\refErrAprOne}{\ensuremath{3.4{\times}10^{-11}}}
\newcommand{\parN}{36{,}481}
\newcommand{\parPmax}{16}
\newcommand{\parGTmax}{7.66}
\newcommand{\parGErr}{\ensuremath{8.1{\times}10^{-11}}}
\newcommand{\parGStages}{56}
\newcommand{\parAdPTone}{15.65}
\newcommand{\parAdPTmax}{2.48}
\newcommand{\parAdPErr}{\ensuremath{5.5{\times}10^{-12}}}
\newcommand{\parAdPStages}{28}
\newcommand{\parAdPSetup}{1.22}
\newcommand{\parAdPLoop}{1.24}
\newcommand{\parHPTone}{23.82}
\newcommand{\parHPTmax}{3.05}
\newcommand{\parHPErr}{\ensuremath{2.5{\times}10^{-13}}}
\newcommand{\parHPStages}{66}
\newcommand{\parHPSetup}{1.97}
\newcommand{\parHPLoop}{1.08}
\newcommand{\parRounds}{6}
\newcommand{\parCondS}{\ensuremath{6.9{\times}10^{11}}}
\newcommand{\parSpread}{6}
\newcommand{\parSchurT}{35.67}
\newcommand{\parSchurErr}{\ensuremath{3.9{\times}10^{-12}}}
\newcommand{\parAdFirstShifts}{91}
\newcommand{\parAdFirstReal}{63}
\newcommand{\parHPShifts}{8}
\newcommand{\parHPRoundsPer}{3}
\newcommand{\parHPWindows}{22}
\newcommand{\parHPSubLengths}{14}
\newcommand{\parHPReal}{1}
\newcommand{\plateNmin}{2{,}004}
\newcommand{\plateNmax}{61{,}234}
\newcommand{\plateLamLast}{\ensuremath{5.7{\times}10^{5}}}
\newcommand{\plateRefExact}{\ensuremath{1.8{\times}10^{-12}}}
\newcommand{\plateGradedErr}{\ensuremath{9.6{\times}10^{-12}}}
\newcommand{\plateAdPErr}{\ensuremath{9.6{\times}10^{-12}}}
\newcommand{\plateHPcErr}{\ensuremath{3.4{\times}10^{-13}}}
\newcommand{\plateLegFirst}{\ensuremath{3.1{\times}10^{-10}}}
\newcommand{\plateLegLast}{\ensuremath{3.2{\times}10^{-4}}}
\newcommand{\plateAprFirst}{\ensuremath{9.1{\times}10^{-12}}}
\newcommand{\plateAprLast}{\ensuremath{1.5{\times}10^{-10}}}
\newcommand{\plateGradedFirst}{\ensuremath{9.6{\times}10^{-12}}}
\newcommand{\plateFirstWindG}{\ensuremath{2.4{\times}10^{-2}}}
\newcommand{\plateFirstWinRadau}{\ensuremath{2.8{\times}10^{-2}}}
\newcommand{\plateBoundMin}{5.8}
\newcommand{\plateBoundMax}{\ensuremath{4.1{\times}10^{1}}}
\newcommand{\sharpZ}{886}
\newcommand{\sharpTerm}{\ensuremath{7.6{\times}10^{-11}}}
\newcommand{\sharpDiff}{\ensuremath{8{\times}10^{-3}}}
\newcommand{\stiffN}{100{,}565}
\newcommand{\stiffLk}{\ensuremath{4.4{\times}10^{6}}}
\newcommand{\stiffGraded}{\ensuremath{1.9{\times}10^{-11}}}
\newcommand{\stiffLeg}{\ensuremath{5.2{\times}10^{-3}}}
\newcommand{\stiffdG}{\ensuremath{1.1{\times}10^{-1}}}
\newcommand{\waveN}{4{,}418}
\newcommand{\waveelmE}{\ensuremath{8.1{\times}10^{-5}}}
\newcommand{\waveelmerr}{\ensuremath{1.2{\times}10^{-4}}}
\newcommand{\waveelmtime}{7.6}
\newcommand{\waveelmlsE}{\ensuremath{1.5{\times}10^{-6}}}
\newcommand{\waveelmlserr}{\ensuremath{1.7{\times}10^{-6}}}
\newcommand{\waveelmlstime}{17.3}
\newcommand{\waveelmlsgrowth}{\ensuremath{2.5{\times}10^{-7}}}
\newcommand{\wavedgrandE}{\ensuremath{8.9{\times}10^{-11}}}
\newcommand{\wavedgranderr}{\ensuremath{1.1{\times}10^{-10}}}
\newcommand{\wavedgrandtime}{7.9}
\newcommand{\wavedggradederr}{\ensuremath{6.6{\times}10^{-9}}}
\newcommand{\wavedglegerr}{\ensuremath{1.8{\times}10^{-8}}}
\newcommand{\wavedglegTwentySevenerr}{\ensuremath{8.7{\times}10^{-12}}}
\newcommand{\wmSeeds}{20}
\newcommand{\wmMeshB}{96}
\newcommand{\wmOmegaA}{24}
\newcommand{\wmOmegaB}{49}
\newcommand{\wmDGAmed}{\ensuremath{7.2{\times}10^{-11}}}
\newcommand{\wmDGAhi}{\ensuremath{1.3{\times}10^{-9}}}
\newcommand{\wmCollAmed}{\ensuremath{8.7{\times}10^{-6}}}
\newcommand{\wmCollAhi}{\ensuremath{1.4{\times}10^{-4}}}
\newcommand{\wmLSAmed}{\ensuremath{1.0{\times}10^{-7}}}
\newcommand{\wmLSAhi}{\ensuremath{1.7{\times}10^{-6}}}
\newcommand{\wmDGBhi}{\ensuremath{1.3{\times}10^{-9}}}
\newcommand{\wmCollBdiv}{15}
\newcommand{\wmRadBdiv}{17}
\newcommand{\wmLSBmed}{\ensuremath{6.2{\times}10^{-8}}}
\newcommand{\grNeurons}{130}
\newcommand{\grDps}{70}
\newcommand{\grCheckDps}{90}
\newcommand{\grGdef}{\ensuremath{8{\times}10^{-46}}}
\newcommand{\grLeb}{4.0}
\newcommand{\grSeconds}{42}
\newcommand{\dblSymMin}{\ensuremath{2.7{\times}10^{1}}}
\newcommand{\dblNeurons}{72}
\newcommand{\dblSmooth}{\ensuremath{3.8{\times}10^{-9}}}
\newcommand{\hpASmooth}{\ensuremath{9.8{\times}10^{-13}}}
\newcommand{\designA}{\ensuremath{2.5{\times}10^{-9}}}
\newcommand{\designB}{\ensuremath{2.3{\times}10^{-9}}}
\newcommand{\designC}{\ensuremath{4.1{\times}10^{-10}}}
\newcommand{\designD}{\ensuremath{4.7{\times}10^{-12}}}
\newcommand{\tailKtwo}{732}
\newcommand{\tailTen}{\ensuremath{9.4{\times}10^{-5}}}
\newcommand{\tailTwelve}{56}
\newcommand{\tailSixteen}{732}
\newcommand{\designDdim}{141}
\newcommand{\lsCond}{\ensuremath{7{\times}10^{27}}}
\newcommand{\lsMaxR}{1.02}
\newcommand{\hoDeg}{6}
\newcommand{\hoN}{20{,}113}
\newcommand{\hoNhi}{35{,}809}
\newcommand{\hoLk}{\ensuremath{1.1{\times}10^{5}}}
\newcommand{\hoSpaceMax}{\ensuremath{4.6{\times}10^{-8}}}
\newcommand{\hoSpaceMin}{\ensuremath{1.7{\times}10^{-10}}}
\newcommand{\hoGradedMax}{\ensuremath{4.1{\times}10^{-12}}}
\newcommand{\hodGFirst}{\ensuremath{2.6{\times}10^{-3}}}
\newcommand{\hoRadauFirst}{\ensuremath{3.8{\times}10^{-4}}}
\newcommand{\hoAdPMax}{\ensuremath{4.1{\times}10^{-12}}}
\newcommand{\hoHPcMax}{\ensuremath{1.9{\times}10^{-12}}}
\newcommand{\hoAprMax}{\ensuremath{1.3{\times}10^{-11}}}
\newcommand{\adSmax}{10^{3}}
\newcommand{\adR}{40}
\newcommand{\adUnk}{399}
\newcommand{\adModal}{\ensuremath{4.3{\times}10^{-10}}}
\newcommand{\adModalLater}{\ensuremath{6.7{\times}10^{-12}}}
\newcommand{\grModalLater}{\ensuremath{4.9{\times}10^{-11}}}
\newcommand{\hpcQ}{14}
\newcommand{\hpcL}{15}
\newcommand{\hpcUnk}{330}
\newcommand{\hpcModal}{\ensuremath{3.1{\times}10^{-11}}}
\newcommand{\hpcTwelveL}{15}
\newcommand{\hpcTwelveUnk}{286}
\newcommand{\hpcTwelveModal}{\ensuremath{5.4{\times}10^{-11}}}
\newcommand{\hpcTenMin}{\ensuremath{1.8{\times}10^{-9}}}
\newcommand{\aprModal}{\ensuremath{1.3{\times}10^{-8}}}
\newcommand{\adPUnk}{224}
\newcommand{\adPModalLater}{\ensuremath{6.7{\times}10^{-12}}}

\begin{document}

\maketitle

\begin{abstract}
Linearized neural networks (LNNs), shallow networks with a fixed hidden layer, have been used as collocation-based
time integrators for stiff problems. We use the hidden layer instead as
the trial space of discontinuous Galerkin time stepping for finite element discretizations of stiff linear evolution
equations. The method is A- and L-stable for every hidden layer, and each window reduces to independent shifted
problems. For symmetric, coercive, homogeneous problems on uniform windows, its error is bounded, independently of the
mesh, by a worst-case modal error computed on a reference window. The analysis shows that polynomials damp
unresolved stiff modes only weakly, by a factor given by the reproducing kernel.
We therefore propose a $\tanh$ layer whose slopes are graded geometrically at the beginning of the window and whose
matrices in time are formed once in extended precision. On uniform windows, this layer resolves stiff initial layers of
the heat equation to about $10^{-11}$ on all meshes tested, including unstructured, locally refined ones, and
with high-order elements, where polynomial layers fail. Replacing the geometrically refined first window of $hp$ time
stepping with polynomials of degree $\hpcQ$ by a single window with the graded layer reduces the unknowns in time and
sequential rounds of solves. Since it also saves factorizations, this variant is faster on one to
$\parPmax$ processes, although $hp$ time stepping remains more accurate. In a long-time wave test, the Galerkin condition makes
extreme learning machines with the same random hidden layer three to five orders of magnitude more accurate.
\end{abstract}



\section{Introduction}\label{sec:intro}

Time discretizations of evolution equations are mostly based on time marching, e.g., implicit Runge--Kutta methods
such as the Radau IIA methods for stiff problems \cite{HairerWanner1996}. The discontinuous Galerkin (DG) time stepping
method dG($q$) \cite{LesaintRaviart1974,Thomee2006} is A- and L-stable, and it is equivalent, at the nodes and for
linear homogeneous problems with constant coefficients, to the Radau IIA method with $q+1$ stages. Its $hp$-version with
geometric refinement towards $t=0$ converges exponentially for parabolic problems with incompatible initial data
\cite{SchoetzauSchwab2000}, and an inf--sup condition and a nearly best approximation property of dG($q$) were
established in \cite{Saito2021,Gomez2026}. Each time step of dG($q$) reduces to $q+1$ independent complex shifted
problems \cite{Butcher1976,SchoetzauSchwab2000}, where the condition number of the transformation increases rapidly
with $q$ \cite{vonPetersdorffSchwab2004}. Non-polynomial trial spaces in time have been used in functionally fitted
Runge--Kutta methods \cite{Ozawa2001,HoangSidje2008} and continuous Galerkin methods \cite{LiWu2016}, and for linear
problems with time-independent operators, Laplace transform methods also reduce the time discretization to independent
shifted problems \cite{SheenSloanThomee2000,WeidemanTrefethen2007}.

On the other hand, physics-informed neural networks (PINNs) \cite{RaissiPerdikarisKarniadakis2019} fit a deep network
to a differential equation by minimizing its residual over all parameters, either on the space-time domain or combined with implicit Runge--Kutta
methods in time, and they are often trained by marching in time \cite{WangSankaranPerdikaris2024}. Shallow networks
whose hidden layer is fixed in advance and whose output layer is computed from a linear problem are known as extreme
learning machines (ELMs) \cite{HuangZhuSiew2006} and random feature methods \cite{ChenELuo2023}. With deterministic,
quasi-uniformly distributed hidden parameters, they achieve optimal approximation rates and are called linearized
neural networks (LNNs) \cite{MaoXuXu2026,LiMaoXu2026}. In this paper, we use the term LNN for any shallow neural network
whose hidden layer is fixed in advance. LNNs have been applied on space-time domains \cite{ChenELuo2023}, with block
time marching \cite{DongLi2021}, with DG formulations \cite{SunDongWang2022}, and as time integrators for stiff
ordinary differential equations (ODEs) and semi-discretized PDEs
\cite{DeFlorioSchiassiFurfaro2022,FabianiGalarisRussoSiettos2023}, where the output layer is computed from the residual
at collocation points. Note that, for LNNs in space, Galerkin formulations lead to severely ill-conditioned linear
systems, and collocation with least squares solvers is preferred \cite{MaoXuXu2026}.

Unlike finite element methods, these methods do not build the boundary and initial conditions into the trial space,
but add them to the loss as penalty terms, impose them weakly by Nitsche's method \cite{LiaoMing2021}, or satisfy them
by an ansatz \cite{LagarisLikasFotiadis1998,LuEtAl2021}, e.g., the initial value of an ODE by the theory of functional connections
\cite{Mortari2017,SchiassiEtAl2021,DeFlorioSchiassiFurfaro2022}. For Dirichlet conditions, such an ansatz requires an
approximate distance function to the boundary \cite{SukumarSrivastava2022}, whose construction is a modeling step of its
own for complex geometries with mixed boundary conditions. Moreover, the stability of collocation with fixed basis
functions depends on the basis and on the collocation points \cite{HoangSidje2008,FabianiBolltSiettosYannacopoulos2024}.

Hence, a natural approach is to discretize the space variable by the finite element method (FEM), which imposes the
boundary conditions exactly on unstructured meshes, and to use an LNN only for the time variable, as in
\cite{FabianiGalarisRussoSiettos2023} for a finite difference discretization. If $u_h(t)=\sum_{j=1}^N u_j(t)\varphi_j$ is
a finite element function, an LNN $t\mapsto\R^N$ whose outputs are its coefficients turns the PDE into the
system of ODEs
\begin{equation}\label{eq:ode}
  M u'(t)+K u(t)=b(t),\qquad u(0)=u_0,
\end{equation}
where $M$ and $K$ are the mass and stiffness matrices, respectively. However, there are several issues in this approach.
First, the output layer and the initial value have to be determined such that the time stepping is stable for a general
hidden layer. Second, the accuracy may depend on the hidden layer and on the mesh, since \eqref{eq:ode} becomes stiffer
as the mesh is refined, and it is not obvious whether the width of the network must grow with $N$. Finally, it is not
clear whether an LNN in time has any advantage over Galerkin time stepping with polynomials.

In this paper, we determine the output layer of an LNN in time by a DG condition, so that the hidden layer is the trial
space of DG time stepping, and we refer to this method as DG-LNN. For polynomial hidden layers, it is dG($q$), and for
other hidden layers, it is a functionally fitted DG method. Testing with the trial space makes DG-LNN A-stable for every
hidden layer, and the upwind term, which imposes the initial value weakly instead of as a hard constraint, adds
L-stability. The analysis also shows that, among all hidden layers of a given dimension, the polynomials attain the
largest order as $\lambda k\to0$, but on uniform windows, they damp the stiff modes that they do not resolve only by a
factor proportional to $(q+1)/(\lambda k)$. Motivated by this, we propose a graded $\tanh$ layer, whose neurons have
their transition at the beginning of the window and slopes graded geometrically over thirteen decades. It moves the
geometric refinement of $hp$-DG from the windows to the parameters of the hidden layer, and its matrices in time have
to be formed once in extended precision, independently of the problem.

Our contributions are summarized as follows:
\begin{itemize}
\item We prove that DG-LNN is A- and L-stable for every hidden layer and uniquely solvable for coercive problems, that
it is consistent if a matrix in time, which depends only on the hidden layer, is nonsingular, and that each window
reduces to independent shifted problems (\cref{sec:algebra,sec:analysis}).
\item For symmetric, coercive, homogeneous problems on uniform windows, we bound the error for all spatial
discretizations by the worst-case modal error of the hidden layer, which is computed on a reference window and can be
used to choose hidden layers in advance, and we show that the reproducing kernel of the hidden layer determines the
damping of unresolved modes (\cref{sec:analysis}).
\item We compare the graded layer with polynomial and other $\tanh$ layers and Runge--Kutta methods on the same uniform
windows and, used only on the first window, with $hp$-DG on one and several processes, and the DG condition with
collocation and least squares for the same hidden layers, for the heat and wave equations (\cref{sec:numerics}). The comparison with $hp$-DG of the same degree $q=\hpcQ$ shows a trade-off: the graded
layer needs fewer factorizations and is faster on one to $\parPmax$ processes, whereas $hp$-DG is more accurate and,
with $q=12$, also faster on one process.
\end{itemize}

The rest of this paper is organized as follows. In \cref{sec:prelim}, we introduce the semi-discrete problem and LNNs in
time. We present DG-LNN and the hidden layers in \cref{sec:method}, and we prove its stability and error bounds in
\cref{sec:analysis}. Numerical results are presented in \cref{sec:numerics}. We conclude this paper with remarks in
\cref{sec:conclusion}.
\section{Preliminaries}\label{sec:prelim}
In this section, we briefly introduce the semi-discrete problem obtained by the FEM and the LNNs in time that are used
throughout this paper.

\subsection{Semi-Discrete Problem}
Let $V_h$ be a finite element space of dimension $N$ equipped with the $L^2$ inner product $\ip\cdot\cdot$ and the norm
$\norm\cdot$, and let $a$ be a time-independent bilinear form on $V_h$ such that
\begin{equation}\label{eq:accretive}
  a(v,v)\ge0\qquad\forall v\in V_h .
\end{equation}
When $a$ is symmetric, we write $\anorm v:=a(v,v)^{1/2}$. For given $u_{0,h}\in V_h$ and $f$, the semi-discrete solution
$u_h\in H^1(0,T;V_h)$ satisfies $u_h(0)=u_{0,h}$ and
\begin{equation}\label{eq:semi}
  \ip{u_h'(t)}{v}+a(u_h(t),v)=\ip{f(t)}{v}\qquad\forall v\in V_h,\ t\in(0,T].
\end{equation}
In a nodal basis, \eqref{eq:semi} is written as \eqref{eq:ode}, where $M$ is symmetric positive definite and $K+K^\top\succeq0$.
Note that \eqref{eq:accretive} covers the heat equation and the wave equation as a first-order system with the energy
inner product (\cref{sec:ex4}).

\subsection{Linearized Neural Networks in Time}
Let $0=t_0<t_1<\dots<t_W=T$ be a partition of $[0,T]$ into $W$ windows $I_n=(t_{n-1},t_n]$ of length
$k_n=t_n-t_{n-1}$, and let $F_n$ be the increasing affine map from the reference window $\hI=(-1,1]$ onto $I_n$.

\begin{definition}[hidden layer]\label{def:layer}
For given activation functions $\sigma_1,\dots,\sigma_m$ and fixed parameters $\{(w_i,b_i)\}_{i=1}^m$, the hidden layer
is the space
\[
  \hS:=\spanop\{1,\ \sigma_1(w_1\tau+b_1),\dots,\sigma_m(w_m\tau+b_m)\}\subset H^1(\hI),\qquad r+1:=\dim\hS,
\]
and we set $S_n:=\{s\circ F_n^{-1}:s\in\hS\}$. Then a function $U\in S_n\otimes V_h$ of the form
\[
  U(t)=\sum_{j=0}^r\phi_j\big(F_n^{-1}(t)\big)\,x_j,\qquad x_j\in V_h\cong\R^N,
\]
where $\{\phi_j\}$ is a basis of $\hS$, is a shallow neural network $t\mapsto\R^N$ with the hidden layer $\hS$ and the
output weights $X=(x_0,\dots,x_r)$. We call it an LNN since the hidden parameters are not optimized, so that $U$ is linear
in $X$.
\end{definition}

Note that the analysis in this paper uses only that $\hS$ is a finite-dimensional subspace of $H^1(\hI)$ that contains
the constant functions. In particular, the polynomial neurons $\sigma_i(\tau)=\tau^i$ with $w_i=1$ and $b_i=0$ give the trial space $\hS=\Pq_q$, $q=m$, of the classical methods. More generally, every finite-dimensional trial
space
that contains the constants is an LNN in this sense, and the term refers to its parametrization by neurons. The hidden layer may also differ between the windows, in which case $S_n$ is defined by a hidden layer $\hS_n$, e.g., $\hS_n=\Pq_{q_n}$ in
$hp$-DG \cite{SchoetzauSchwab2000}. Unless stated otherwise, we use the same hidden layer on every window.
\section{Galerkin Linearized Neural Networks in Time}\label{sec:method}
In this section, we present DG-LNN, which determines the output layer of an LNN in time by a DG condition, derive its
algebraic structure, and introduce the hidden layers used in this paper.

\subsection{Discontinuous Galerkin Formulation}\label{sec:algebra}
Let $U(t_0^-):=u_{0,h}$. For given $U(t_{n-1}^-)$, DG-LNN finds $U|_{I_n}\in S_n\otimes V_h$ such that
\begin{equation}\label{eq:dg}
\begin{aligned}
  B_n(U,V)&:=\int_{I_n}\!\big[\ip{U'}{V}+a(U,V)\big]\dd t+\ip{U(t_{n-1}^+)}{V(t_{n-1}^+)}\\
  &\phantom{:}=\int_{I_n}\!\ip fV\dd t+\ip{U(t_{n-1}^-)}{V(t_{n-1}^+)}\qquad\forall V\in S_n\otimes V_h .
\end{aligned}
\end{equation}
For $\hS=\Pq_q$, \eqref{eq:dg} is the DG time stepping method dG($q$) \cite{LesaintRaviart1974,Thomee2006}. The initial
value enters through the upwind term and is carried by the constant function in $\hS$. In contrast to a hard constraint, it is imposed weakly, so that
$U(t_{n-1}^+)\ne U(t_{n-1}^-)$ in general, and the jump absorbs
the components that the hidden layer does not resolve (\cref{sec:analysis}). Note that the ELM time
integrators \cite{DeFlorioSchiassiFurfaro2022,FabianiGalarisRussoSiettos2023,DongLi2021} use a random hidden layer and
determine the output layer from the initial value and the residual of \eqref{eq:ode} at collocation points in each
window instead of \eqref{eq:dg}.

Let $\{\phi_j\}_{j=0}^r$ be a basis of $\hS$, let $\phi(\tau)$ be the vector of its values, and define
\begin{equation}\label{eq:GA}
  G_{ij}=\int_{\hI}\phi_i\phi_j\dd\tau,\quad A_{ij}=\int_{\hI}\phi_i\phi_j'\dd\tau+\phi_i(-1)\phi_j(-1),\quad
  D_n:=\frac{2}{k_n}G^{-1}A .
\end{equation}
Let $X\in\R^{(r+1)\times N}$ have the rows $x_j^\top$. Testing \eqref{eq:dg} with $V=\phi_i(F_n^{-1}\cdot)\,v$,
$v\in V_h$, and multiplying by $\tfrac2{k_n}G^{-1}$, we obtain the Sylvester equation
\begin{equation}\label{eq:sylvester}
  D_n X M^\top+X K^\top=Y_n,
\end{equation}
where $Y_n:=\tfrac2{k_n}G^{-1}\phi(-1)\,(MU(t_{n-1}^-))^\top+G^{-1}\int_{\hI}\phi(\tau)\,b(F_n(\tau))^\top\dd\tau$.

As for implicit Runge--Kutta methods \cite{Butcher1976} and dG($q$) \cite{SchoetzauSchwab2000}, \eqref{eq:sylvester}
is decoupled by a diagonalization of the matrix in time. Assume that $D_n$ is diagonalizable, $D_n=\mathcal S\diag(\mu_0,\dots,\mu_r)\mathcal S^{-1}$ (otherwise, the Schur form below is used). Multiplying \eqref{eq:sylvester} by $\mathcal S^{-1}$ from the
left, we obtain the $r+1$ independent problems
\begin{equation}\label{eq:shifted}
  (\mu_jM+K)\,\tilde x_j=\tilde y_j,\qquad j=0,\dots,r,\qquad \tilde X=\mathcal S^{-1}X,\quad \tilde Y=\mathcal S^{-1}Y_n,
\end{equation}
where $\tilde x_j^\top$ and $\tilde y_j^\top$ are the rows of $\tilde X$ and $\tilde Y$, respectively. For real data, $\mu_k=\bar\mu_j$
implies $\tilde x_k=\overline{\tilde x_j}$, so that one solve is required for each real eigenvalue and for each pair of complex conjugate
eigenvalues. Furthermore, the matrices
$\mu_jM+K$ depend only on $(\hS,k_n)$, so that their factorizations (or preconditioner setups) are reused in every
window of the same length.

We show in \cref{thm:spectral} that $\Re\mu_j\ge0$ for every hidden layer, so that all shifted matrices in
\eqref{eq:shifted} are nonsingular when $a$ is coercive. In floating point arithmetic, however, the rounding errors of
the shifted solves are amplified by the condition number $\kappa(\mathcal S)$ of the transformation, which increases
rapidly with the degree for polynomial hidden layers \cite{vonPetersdorffSchwab2004}. Hence, we apply iterative
refinement to the residual of \eqref{eq:sylvester}, where each step requires one more set of the independent
back-solves but no new factorization. As a rule of thumb, we use
$\ell_{\rm ref}=\lceil\log(10^{-14})/\log(10\,\kappa(\mathcal S)\varepsilon_{\rm solve})\rceil-1$ steps, provided that
$10\,\kappa(\mathcal S)\varepsilon_{\rm solve}<1$, with the nominal relative accuracy $\varepsilon_{\rm solve}=10^{-15}$
of the shifted solves by sparse LU. Alternatively, the complex Schur form $D_n=QTQ^{\mathsf H}$ gives the sequential
Bartels--Stewart algorithm \cite{BartelsStewart1972}, as in \cite{SchoetzauSchwab2000,vonPetersdorffSchwab2004} for
dG($q$), with $r+1$ sequential solves per window. Since $\norm{D_n}$ is large for the graded layer
(\cref{rem:implementation}), we apply one step of iterative refinement after it as well. For example, for the heat
equation in \cref{sec:ex3} on a $64^2$ mesh with the graded layer and $W=8$ windows, the maximal relative error is
$\refErrSchurZero$ with the Schur form, $\refErrSchurOne$ after one step of refinement, and $\refErrSchurTwo$ after two
steps, whereas with the diagonalization ($\kappa(\mathcal S)=\refCondS$), it is $\refErrDiagFour$, $\refErrDiagSix$, and
$\refErrDiagEight$ after $4$, $6$, and $8$ steps, respectively. The errors in \cref{sec:ex3,sec:plate} show that one
step suffices up to $N=\stiffN$, and for $hp$-DG with the a priori plan in \cref{sec:setup}, one step after the Schur
form changes the error from $\refErrAprZero$ to $\refErrAprOne$. Hence, in serial runs, we use the Schur form with one
step of refinement for all DG methods whose $\ell_{\rm ref}$ exceeds one, and the diagonalization (\cref{alg:window}) for
parallel runs.

\begin{algorithm}[t]
\caption{DG-LNN time stepping}\label{alg:window}
\begin{algorithmic}
\State \textbf{Setup} (once per window length $k$): form $G,A$ by quadrature, with $\operatorname{sym}(A)$ set to its exact value
\State \quad $\tfrac12(\phi(1)\phi(1)^\top+\phi(-1)\phi(-1)^\top)$, and compute $D=\tfrac2kG^{-1}A=\mathcal S\diag(\mu_j)\mathcal S^{-1}$
\State \quad factorize the complex matrices $\mu_jM+K$ for one $\mu_j$ of each conjugate pair
\For{$n=1,\dots,W$}
\State assemble $Y_n$ from $U(t_{n-1}^-)$ and $f$, and set $X\gets0$
\For{$\ell=0,\dots,\ell_{\mathrm{ref}}$}
\State $\tilde Y\gets\mathcal S^{-1}(Y_n-DXM^\top-XK^\top)$
\State solve $(\mu_jM+K)\tilde x_j=\tilde y_j$ for all $j$ \textit{in parallel}, and $X\gets X+\mathcal S\tilde X$
\EndFor
\State $U|_{I_n}\gets\sum_j\phi_j(F_n^{-1}\cdot)x_j$ and $U(t_n^-)\gets X^\top\phi(1)$
\EndFor
\end{algorithmic}
\end{algorithm}

\subsection{Hidden Layers}\label{sec:layers}
We compare the following hidden layers:
\begin{itemize}
\item \emph{Legendre}: $\hS=\Pq_q$, which gives the classical dG($q$); we use $q=16$ and, for comparisons at the size of
  the graded layer, $q=\univr$.
\item \emph{tanh, quasi-uniform} \cite{MaoXuXu2026}: $(w_i,b_i)=\rho(\cos\theta_i,\sin\theta_i)$ with equally spaced
  $\theta_i$ on a half circle and $\rho=8$; we use $m=32$ and $m=200$ neurons, for which the compressed dimension
  saturates at $r+1=83$.
\item \emph{tanh, random}: $(w_i,b_i)\sim U[-4,4]^2$ with a fixed seed, as in ELMs
  \cite{HuangZhuSiew2006,DongLi2021}; we use $m=32$ ($r+1=28$ after compression) and
  $m=192$, for which the compressed dimension saturates at $r+1=46$.
\item \emph{tanh, graded}: the neurons \eqref{eq:graded} below.
\end{itemize}
The radius $\rho$ and the range of the random parameters are scale parameters, and with our values, the slopes are at
most about $8$, so that these layers suit smooth functions rather than stiff decay.
The graded $\tanh$ layer consists of the neurons
\begin{equation}\label{eq:graded}
  \tanh\big(s_i(\tau+1)\big),\qquad s_i=s_{\min}(s_{\max}/s_{\min})^{(i-1)/(m-1)},\quad i=1,\dots,m,
\end{equation}
with $[s_{\min},s_{\max}]=[10^{-2},\,10^{11}]$ and $m=\lceil\nu\log_{10}(s_{\max}/s_{\min})\rceil$, so that the slopes are geometrically spaced with $\nu$ slopes per decade. We use $\nu=10$, which gives $m=\grNeurons$. This is a pure LNN, which contains no polynomials except the constants, and its parameters $(w_i,b_i)$ lie on the ray $w=b$ with geometrically graded radii. In contrast to the quasi-uniform parameters of \cite{MaoXuXu2026,LiMaoXu2026}, all
neurons have their transition at the beginning of the window, and their widths $1/s_i$ are graded geometrically, which
is the counterpart of the geometric refinement of the windows towards $t=0$ in $hp$-DG \cite{SchoetzauSchwab2000}. The
motivation is the expansion $\tanh(s(\tau+1))=1-2e^{-2s(\tau+1)}+2e^{-4s(\tau+1)}-\cdots$, so that the neuron with slope $s$ represents the decaying exponential $e^{-\lambda(t-t_{n-1})}$ with $\lambda k\approx 4s$, which is the time
dependence of a stiff mode of \eqref{eq:semi}. For pure exponentials, the relative distance from $e^{-\lambda t}$ to
$\spanop\{e^{-s_1t},\dots,e^{-s_nt}\}$ in $L^2(0,\infty)$ is $\prod_{i=1}^n|\lambda-s_i|/(\lambda+s_i)$ by the Cauchy
determinant, and minimizing it uniformly over $\lambda\in[a,b]$ is a Zolotarev problem, whose value is bounded by a
multiple of $\exp(-\pi^2n/(2\log(4b/a)))$ \cite[Appendix~A]{BeckermannTownsend2017}. Since the product is invariant
under a common scaling of $\lambda$ and the rates, we space the slopes geometrically. Relaxation times sampled from a
log-uniform law were also used for the temporal features in \cite{LinghuWang2026}. Note that the neurons are not pure
exponentials, so that this is a motivation and not an error estimate. Since the slopes cover
$\lambda k\in[0.04,\,4\cdot10^{11}]$, the graded layer is fixed in advance for this range, without knowledge of the
operator or the mesh, and modes with larger $\lambda k$ are damped (\cref{rem:tail}). It does not depend on the window
length, but the window length has to resolve the forcing and the smooth part of the solution (\cref{sec:ex3}). Our
choices of $\nu$, the slope range, and the truncation tolerance below are discussed in \cref{sec:ex1}.

\begin{remark}[implementation]\label[remark]{rem:implementation}
The neurons \eqref{eq:graded} are nearly linearly dependent, and their orthonormalization involves coefficients of the
size of the inverse of the smallest singular value, which amplify the rounding errors in double precision, so that
the computed $G$ and $A$ violate \cref{thm:spectral}(i). For example, for the
$\dblNeurons$ neurons with $\nu=8$ and $s_{\max}=10^7$, compressed and integrated in double precision as the other
hidden layers, the symmetric part of $A$ has a negative eigenvalue of modulus $\dblSymMin$, and after the replacement of
the symmetric part by its exact value, the error of a smooth mode ($\lambda=0.1\pi^2$, $T=2$, one or eight windows) inside
the windows is $\dblSmooth$, compared with $\hpASmooth$ with the construction below. For the neurons \eqref{eq:graded},
the double-precision construction breaks down. Hence, we form $G$ and $A$ for $\{1,\tanh(s_i(\tau+1))\}$ with
$\grDps$ digits by a composite Gauss--Legendre rule on a partition graded geometrically
towards $\tau=-1$, and we orthonormalize their span by an eigendecomposition of $G$ with the same precision, truncated at the relative tolerance $10^{-32}$ for the eigenvalues ($10^{-16}$ for the singular values). Then $\hS$ has the
dimension $r+1=\subUnknowns$, $G=I$ and the symmetric part of $A$ is exact up to $\grGdef$, and a finer rule with
$\grCheckDps$ digits does not change $A$ in double precision. Only $A$, $\phi(\pm1)$, and the values of $\phi$ at the
points that are needed are rounded to double precision, and the symmetric part of the rounded $A$ is set to its exact
value, so that the proof of \cref{thm:spectral} applies to the computed matrices. This construction is done once,
independently of the problem, and takes $\grSeconds$\,s. The data $f$ are replaced by their interpolant in $\hS$ at
$r+1$ points selected by a column-pivoted QR factorization, whose Lebesgue constant is at most $\grLeb$, so that
\eqref{eq:energy} holds for this modified problem. For the other hidden layers, we orthonormalize their span by a
truncated singular value decomposition with relative tolerance $10^{-12}$ in double precision, since their symmetric
parts are exact up to rounding errors. Note that $\norm{A}\approx\gradedNormA$ for the graded layer, which amplifies
the rounding errors of the shifted solves (\cref{sec:ex1}).
\end{remark}
\section{Stability and Error Analysis}\label{sec:analysis}

In this section, we prove that DG-LNN is A- and L-stable for every hidden layer, and we bound its error on uniform
windows for all spatial discretizations by a quantity computed on the reference window. Throughout this section, $\hS$
is an arbitrary hidden layer in the sense of \cref{def:layer}, and no property of the activation functions or
hidden parameters is used. For the scalar problem $y'=\lambda y$ with $\Re\lambda\le0$ ($V_h=\C$, $a(u,v)=-\lambda u\bar v$, and $f=0$, where
\eqref{eq:accretive} is understood for the real part), the DG-LNN solution of
one window of length $k$ satisfies $U(t_1^-)=R_S(\zeta)\,y(0)$ with $\zeta=\lambda k$. We call $R_S$ the
\emph{amplification factor}.

We first recall the energy identity of DG time stepping. Let \eqref{eq:accretive} hold, let $U$ solve \eqref{eq:dg}, and
let $[U]_j:=U(t_j^+)-U(t_j^-)$. Taking $V=U|_{I_n}$ in \eqref{eq:dg}, we obtain
\begin{equation}\label{eq:energy}
  \tfrac12\norm{U(t_n^-)}^2+\tfrac12\norm{[U]_{n-1}}^2+\int_{I_n}a(U,U)\dd t=\tfrac12\norm{U(t_{n-1}^-)}^2+\int_{I_n}\ip fU\dd t
\end{equation}
exactly as for dG($q$) \cite[Chap.~12]{Thomee2006}, since this computation requires only that $U|_{I_n}$
is an admissible test function. Hence, \eqref{eq:energy} holds for every hidden layer. In particular,
$\norm{U(t_n^-)}\le\norm{u_{0,h}}$ for all $n$ if $f=0$.

\begin{theorem}[spectrum, solvability, A- and L-stability]\label{thm:spectral}
Let $G,A,D_n$ be as in \eqref{eq:GA}. For every hidden layer:
\begin{enumerate}[(i)]
\item $A+A^\top=\phi(1)\phi(1)^\top+\phi(-1)\phi(-1)^\top\succeq0$, hence $\sigma(D_n)\subset\{\mu\in\C:\Re\mu\ge0\}$.
\item The window problem \eqref{eq:dg} has a unique solution if $a$ is coercive on $V_h$, or if \eqref{eq:accretive} holds and
  $\sigma(D_n)\cap i\R=\emptyset$. The latter is a property of $\hS$ alone. In both cases all shifted matrices
  $\mu_jM+K$ in \eqref{eq:shifted} are nonsingular.
\item $R_S(\zeta)=\phi(1)^\top\big(A-\tfrac\zeta2G\big)^{-1}\phi(-1)$ is a rational function whose numerator and
  denominator have degrees at most $r$ and $r+1$, respectively. It satisfies $|R_S(\zeta)|\le1$ for $\Re\zeta\le0$
  (A-stability) and $R_S(\zeta)=O(\zeta^{-1})$ as $|\zeta|\to\infty$ (L-stability).
\end{enumerate}
\end{theorem}
\begin{proof}
(i) Integration by parts gives $\int_{\hI}(\phi_i\phi_j'+\phi_i'\phi_j)=\phi_i(1)\phi_j(1)-\phi_i(-1)\phi_j(-1)$, which yields
the formula for $A+A^\top$. If $Av=\mu'Gv$ with $v\in\C^{r+1}\setminus\{0\}$, then
$\Re\mu'\,v^HGv=\Re v^HAv=\tfrac12(|\phi(1)^\top v|^2+|\phi(-1)^\top v|^2)\ge0$. Since $G\succ0$, $\Re\mu'\ge0$, and
$\sigma(D_n)=\tfrac2{k_n}\sigma(G^{-1}A)$.

(ii) It is well known that \eqref{eq:sylvester} is uniquely solvable if and only if $\mu M+K$ is nonsingular for every
$\mu\in\sigma(D_n)$. For $w\in\C^N\setminus\{0\}$,
$\Re w^H(\mu M+K)w=\Re\mu\,w^HMw+w^HK_sw$ with $K_s=\tfrac12(K+K^\top)\succeq0$. This is positive if $\Re\mu>0$, and also if
$\Re\mu=0$ and $K_s\succ0$ (coercivity).

(iii) In the scalar case, \eqref{eq:dg} reads $(A-\tfrac\zeta2G)x=\phi(-1)y(0)$ with $U(t_1^-)=\phi(1)^\top x$, which gives
the formula, and Cramer's rule gives the degrees, since $\det(A-\tfrac\zeta2G)$ has the leading coefficient
$\det(-\tfrac12G)\ne0$. For $\Re\zeta<0$, the matrix $A-\tfrac\zeta2G$ is nonsingular by (i). Taking the real part of
$x^H(A-\tfrac\zeta2G)x=x^H\phi(-1)\,y(0)$ and using (i), we obtain with $u_\pm:=\phi(\pm1)^\top x$
\[
  |u_+|^2+|u_--y(0)|^2-\Re\zeta\,x^HGx=|y(0)|^2,
\]
so that $|R_S(\zeta)|\le1$. Hence, the rational function $R_S$ has no pole on the imaginary axis, and the bound holds for
$\Re\zeta\le0$ by continuity, where $R_S$ denotes the continuous extension at the points of the imaginary axis at which
$A-\tfrac\zeta2G$ is singular. Finally, $R_S(\zeta)=O(\zeta^{-1})$ follows from the degrees.
\end{proof}

Note that the proofs of (i) and (iii) use only the symmetric part of $A$ and $G\succ0$, so that they also apply to the
computed matrices of \cref{rem:implementation}.

Note that, for $\hS=\Pq_q$, the identity in the proof of (i) was used in \cite{vonPetersdorffSchwab2004} to bound the
eigenvalues of $D_n$. For a general hidden layer, $D_n$ can have eigenvalues on the imaginary axis, e.g.,
$\sigma(G^{-1}A)=\{\tfrac12,\pm i\pi/\sqrt2\}$ for $\hS=\spanop\{1,\,\sin\pi\tau,\,\cos\pi\tau+\cos2\pi\tau\}$, so that
the conditions in (ii) cannot be dropped. Moreover, \cref{thm:spectral} does not imply consistency. Testing
\eqref{eq:dg} for the scalar problem with $V=1$ gives $R_S(\zeta)=1+\tfrac\zeta2\int_{\hI}U\dd\tau$ for the solution $U$
with $y(0)=1$. If $A$ is nonsingular, then $U=1$ for $\zeta=0$, since the constant function solves this problem, so that $R_S(\zeta)=e^\zeta+O(\zeta^2)$ and the order is at least one. Otherwise, this can fail. For example, $A$ is singular
for $\hS=\spanop\{1,\tau^2\}$, and $R_S(\zeta)=1/(1-\zeta/6)$, so that DG-LNN with this hidden layer is A- and L-stable
but not consistent. For the hidden layers of \cref{fig:stability}, $A$ is nonsingular, since $\Re\mu_j\ge\dgMinReMin$ for
$k_n=1$.

Next, we consider symmetric coercive problems with $f=0$ on uniform windows of length $k$. For $z>0$, let
$y_z(\tau):=e^{-z(1+\tau)/2}$, which is the mode $e^{-\lambda(t-t_{n-1})}$ on a window with $z=\lambda k$, and let
$U_z\in\hS$ be the DG-LNN solution on the reference window of the scalar problem whose exact solution is $y_z$, so that
$U_z(1)=R_S(-z)$.

\begin{theorem}[mesh-independent worst-case modal error]\label{thm:modal}
Let $a$ be symmetric and coercive, $f=0$, and let the windows be uniform with length $k$. Then for every finite element
space $V_h$ (any mesh, any $N$, any spectrum of $a$), every $n\ge1$, and every $t=F_n(\tau)$ with $\tau\in(-1,1]$,
\[
  \norm{u_h(t)-U(t)}\le\mathcal E_S(n,\tau)\,\norm{u_{0,h}},
\]
where
\[
  \mathcal E_S(n,\tau):=\sup_{z>0}\big|R_S(-z)^{n-1}U_z(\tau)-e^{-(n-1)z}y_z(\tau)\big|.
\]
In particular, $\mathcal E_S(n):=\mathcal E_S(n,1)=\sup_{z>0}|R_S(-z)^n-e^{-nz}|$ bounds the error at the nodes. If
$u_{0,h}$ is an eigenfunction of $a$ with eigenvalue $\lambda$, the error equals the term in the supremum at $z=\lambda k$
times $\norm{u_{0,h}}$, so that the bound is attained if the supremum is attained at such a $z$.
\end{theorem}
\begin{proof}
Let $\{(\lambda_\ell,w_\ell)\}_{\ell=1}^N$ be the $\ip\cdot\cdot$-orthonormal eigenpairs of $a$ on $V_h$, where
$\lambda_\ell>0$, and let $a_\ell:=\ip{u_{0,h}}{w_\ell}$ and $z_\ell:=\lambda_\ell k$. Since
$S_n\otimes V_h=\bigoplus_\ell S_n\otimes\spanop\{w_\ell\}$ and both forms in \eqref{eq:dg} are diagonal in
$\{w_\ell\}$, testing \eqref{eq:dg} with $V=s\,w_\ell$ shows that $U=\sum_\ell U_\ell w_\ell$, where $U_\ell$ is the
scalar DG-LNN solution of $y'+\lambda_\ell y=0$, $y(0)=a_\ell$. Hence,
$U_\ell(F_n(\tau))=a_\ell R_S(-z_\ell)^{n-1}U_{z_\ell}(\tau)$, whereas
$\ip{u_h(F_n(\tau))}{w_\ell}=a_\ell e^{-(n-1)z_\ell}y_{z_\ell}(\tau)$. Therefore,
\[
  \norm{u_h(t)-U(t)}^2=\sum_\ell\big|R_S(-z_\ell)^{n-1}U_{z_\ell}(\tau)-e^{-(n-1)z_\ell}y_{z_\ell}(\tau)\big|^2|a_\ell|^2,
\]
which gives the bound, with equality for $u_{0,h}=w_\ell$ if the supremum is attained at $z_\ell$. This is the classical spectral
argument for rational one-step methods \cite[Chap.~7]{Thomee2006}.
\end{proof}

Note that $\mathcal E_S(n,\tau)$ is computed on $\hI$ alone (\cref{sec:certificates}), so that the width of the
network does not have to grow with $N$. This independence of $N$ holds for every rational one-step method, and what
depends on the hidden layer is the size of $\mathcal E_S$. The proof also applies if the hidden layer $\hS_n$ and the
length $k_n$ depend on $n$, e.g., for $hp$-DG. Then the term in the supremum is replaced by
$\prod_{m<n}R_{S_m}(-\lambda k_m)\,U^{(n)}_{\lambda k_n}(\tau)-e^{-\lambda t_{n-1}}y_{\lambda k_n}(\tau)$, where $U^{(n)}_z$
is defined with $\hS_n$, and the supremum is taken over $\lambda>0$. We use this bound to choose such hidden layers in
advance in \cref{sec:setup}. Since its supremum is taken on a grid in $z\in[10^{-3},10^{12}]$,
the computed values are not rigorous upper bounds. Note also that $\mathcal E_S(1,\tau)$ cannot be small uniformly in $\tau$. In fact, $U_z\to0$ uniformly on $\hI$
as $z\to\infty$, whereas $y_z(-1+2\epsilon/z)=e^{-\epsilon}$ for all $z$. Hence, the accuracy inside the first window is
meaningful only for $t\ge\rho k$ with a fixed $\rho>0$, e.g., $\rho=10^{-3}$ in \cref{sec:ex1}.

The size of $\mathcal E_S(1)$ is determined by two regimes of $z$. For small $z$, \cref{thm:spectral}(iii)
implies that $R_S(-z)-e^{-z}=O(z^{p+1})$ with $p\le2r+1$, and $p=2r+1$ holds only for the Pad\'e approximant of type
$(r,r+1)$ \cite[Sect.~IV.3]{HairerWanner1996}, which is the amplification factor of dG($r$) \cite{LesaintRaviart1974}. Hence, at a given dimension of the hidden layer, the polynomials attain the largest order as $z\to0$, and a
non-polynomial hidden layer can only trade this order for the approximation of stiff modes. For large $z$, the
modes that are not resolved by the hidden layer are only damped, as shown in the following remark.

\begin{remark}[damping of unresolved modes]\label[remark]{rem:tail}
Let $\phi$ be a basis of $\hS$. Since $G$ is nonsingular, $R_S(-z)=\phi(1)^\top(A+\tfrac z2G)^{-1}\phi(-1)$ satisfies
\[
  R_S(-z)=\frac2z\,K_S(1,-1)+O(z^{-2})\quad\text{as }z\to\infty,\qquad K_S(\sigma,\tau):=\phi(\sigma)^\top G^{-1}\phi(\tau),
\]
where $K_S$ is the reproducing kernel of $\hS$ in $L^2(\hI)$. Hence, the modes that are not resolved by the hidden layer
are damped by the factor $2|K_S(1,-1)|/z$, which is determined by the correlation of the hidden layer between the two
ends of the reference window. For $\hS=\Pq_q$, the Legendre expansion of $K_S$ gives $2K_S(1,-1)=(-1)^q(q+1)$, so that a
higher degree damps a fixed unresolved mode less. For the Legendre, quasi-uniform, and random $\tanh$ layers of
\cref{fig:stability}, $z\,R_S(-z)$ at $z=10^{12}$ agrees with $2K_S(1,-1)$ up to a relative difference of
$\tailOthersDev$. For the graded layer, $|2K_S(1,-1)|=\tailKtwo$, and this asymptotic regime starts only beyond its largest slope, with $|z\,R_S(-z)|=\tailTen$, $\tailTwelve$, and $\tailSixteen$ at $z=10^{10}$, $10^{12}$, and
$10^{16}$, respectively. A hidden layer that approximates $y_z$ up to large $z$, and thus resolves the stiff modes inside the window, needs
neurons whose widths match $1/z$, as in
\eqref{eq:graded}.
\end{remark}
\section{Numerical Experiments}\label{sec:numerics}
In this section, we verify the stability and the error bounds in \cref{sec:analysis}, and we compare DG-LNN with time
stepping methods based on polynomials in pairs with the same flexibility. DG-LNN with the graded $\tanh$ layer on every uniform window is compared with dG($q$), the Radau IIA method, and other hidden layers on the same windows, and
DG-LNN with adapted layers, which differ between the windows, is compared with $hp$-DG, which refines the windows
instead. We report the number of unknowns in time, which is the sum of $r+1$ over the windows, and the relative error
$\norm{u_h(t)-U(t)}/\norm{u_h(t)}$ of the time discretization, where $u_h$ is the solution of \eqref{eq:semi}. For the
heat equation, $u_h$ is computed by $hp$-DG with degrees from $4$ to $\refQ$ on $\refL$ geometric levels (ratio
$\tfrac14$) in the first of $8$ windows and $\Pq_{\refQ}$ on the others, which differs from the exact solution by at most
$\meshRefExact$ on the meshes with at most $4000$ unknowns and from the Radau IIA method of order 5 with $4000$ steps by
at most $\heatRefX$. For the wave equation, $u_h$ is the exact modal
solution. For the heat equation in \cref{sec:ex3}, the spatial discretization error on the $128^2$ mesh,
estimated by the solution on the $256^2$ mesh, is $\heatSpaceErrFirst$ at $t=0.01$ and $\heatSpaceErr$ at $t=1$, far
above the time discretization errors.

\subsection{Implementation Details and Baselines}\label{sec:setup}
We use the Firedrake library \cite{FiredrakeUserManual} for the P1 finite element discretizations, where the Dirichlet
degrees of freedom are eliminated. DG-LNN is implemented in Python using SciPy \cite{SciPy2020}, and the shifted systems
\eqref{eq:shifted} are factorized by its sparse direct solver SuperLU. All computations were performed on a laptop with
an Apple M5 Max and 64GB RAM with one thread per process, and the timings in \cref{sec:ex3} are the medians of three
runs. They include the assembly of the right-hand sides, the matrices in time and their decompositions, and all sparse
factorizations, but not the construction of the hidden layers, which is done once and independently of the problem
(\cref{rem:implementation}). Besides the Crank--Nicolson (CN) method and the 3-stage Radau IIA method of order~5
(Radau IIA(5)), we use the following baselines:
\begin{itemize}
\item \emph{Adaptive ESDIRK method}: the L-stable ESDIRK method of order~5 \cite{KennedyCarpenter2003} in PETSc TS
\cite{PETScTS2018} with the sparse direct solver MUMPS instead of SuperLU, which runs compiled code with Python
callbacks.
\item \emph{$hp$-DG} \cite{SchoetzauSchwab2000}: DG-LNN with Legendre layers on uniform windows of length $k$, where the
first window is refined geometrically towards $t=0$ on $L$ levels with ratio $\sigma$, i.e., into sub-windows with the
nodes $k\sigma^{L-1},\dots,k\sigma,k$. We use the constant degree $q=\hpcQ$, which is the largest degree in the
experiments of \cite{SchoetzauSchwab2000}, with $\sigma=\tfrac12$ and $L=\hpcL$, which is chosen below, and the a priori
plan with $\sigma=\tfrac14$ and degrees increasing linearly from $4$ to $\aprQ$ on $\aprL$ levels. For the computing
times, we also use $q\in\{10,12\}$ with $\sigma=\tfrac12$ and $L=\hpcL$, and degrees increasing linearly from $3$ to $q$
on $L$ levels with $\sigma=\tfrac14$ and $(q,L)\in\{(8,4),(12,6),(16,6),(16,10),(20,12)\}$. Every length of the
sub-windows requires its own factorizations.
\item \emph{Adapted layers}: DG-LNN with the graded layer on the first window and a smaller hidden layer on the other
windows, where the stiff modes have decayed, so that, as in $hp$-DG, only the first window resolves the initial layer. In
graded/$\Pq_{\hpcQ}$, the other windows use $\Pq_{\hpcQ}$ as in $hp$-DG, so that the two methods differ only in the
first window. In graded/truncated, they use the graded layer \eqref{eq:graded} truncated at $s_{\max}=\adSmax$ ($\nu=10$), constructed as in
\cref{rem:implementation} with $r+1=\adR$.
With $8$ windows, this gives $\adPUnk$ and $\adUnk$ unknowns in time, respectively, instead of $8(r+1)=\timeGradedUnk$.
\item \emph{ELM} \cite{DeFlorioSchiassiFurfaro2022,DongLi2021}: the random $\tanh$ layer of \cref{sec:layers} with $32$
neurons ($r=27$) or the neurons of the graded layer on fixed windows, where the initial value of each window is imposed
exactly by the ansatz $U(t)=U(t_{n-1}^-)+\sum_{j=1}^r\psi_j(F_n^{-1}(t))\,x_j$ with a basis $\{\psi_j\}$ of
$\{s\in\hS:s(-1)=0\}$. The output weights are determined either by collocation at $r$ points of each window, which are
selected from Chebyshev points and points graded geometrically towards $t_{n-1}$ by a column-pivoted QR factorization,
or by least squares, where the residual of \eqref{eq:ode} is minimized in the norm induced by $M^{-1}$ at $4r$
Chebyshev--Gauss--Lobatto points, where the normal equations are solved by the complex Schur form of a matrix of size
$2r$ with $2r$ shifted factorizations. These are generic ELM time integrators,
whereas the cited methods use other features, adaptive step sizes, or tuned parameters.
\end{itemize}
The truncation and the number of levels are chosen in advance by the worst-case modal error of \cref{thm:modal} for
hidden layers that differ between the windows, taken as the maximum over the $49$ sample times of \cref{sec:ex3} with $8$ windows and over $z=\lambda k\in[10^{-3},4\cdot10^{11}]$, the range of the graded
layer. We take the smallest $s_{\max}\in\{10,10^2,10^3,10^4\}$ for which this error on the windows $n\ge2$ does not
exceed that of the graded layer on all windows ($\grModalLater$), which gives $s_{\max}=\adSmax$ with $\adModalLater$,
whereas $\Pq_{\hpcQ}$ gives $\adPModalLater$ on these
windows. For both variants, the error is then $\adModal$, attained in the first window. For $hp$-DG with $q=\hpcQ$, we
take the smallest $L$ for which the error does not exceed $\adModal$, which gives $L=\hpcL$ with $\hpcModal$. With $q=12$, $L=\hpcTwelveL$ also suffices ($\hpcTwelveModal$ with $\hpcTwelveUnk$
unknowns), whereas with $q=10$, the error is at least $\hpcTenMin$ for all $L\le18$. At this worst-case modal error, $hp$-DG with $q=\hpcQ$ requires $\hpcUnk$ unknowns in time, compared with $\adPUnk$ for graded/$\Pq_{\hpcQ}$. The
a priori plan gives $\aprModal$. The code and data that reproduce all results are provided as supplementary materials.

\subsection{Stability and Error Bounds of Hidden Layers}\label{sec:ex1}
\begin{figure}[tbhp]
\centering\includegraphics[width=\textwidth]{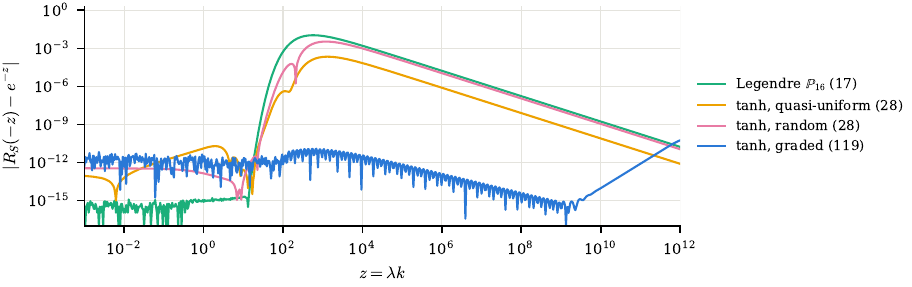}
\caption{Error $|R_S(-z)-e^{-z}|$ of a decaying mode $e^{-\lambda t}$ after one window of DG-LNN, where
$z=\lambda k$. The numbers in the legend are the dimensions of the hidden layers.}\label{fig:stability}
\end{figure}

For the Legendre layer $\Pq_{16}$ ($17$ functions) and the quasi-uniform and random $\tanh$ layers with $32$ neurons
($28$ functions each), DG-LNN satisfies $\max_{\Re\zeta\le0}|R_S|-1\le\stabGalMax$ on our sample of the left half-plane
(\cref{sec:certificates}). For the graded layer ($\subUnknowns$ functions), the computed values exceed one by up to
$\stabGradedMax$, which is the size of the rounding errors amplified by $\norm A\approx\gradedNormA$, whereas the exact
amplification factor of the computed matrices satisfies $|R_S|\le1$ by \cref{thm:spectral}. All four layers satisfy
$|R_S(-10^{12})|\le\dgRinfMax$. \Cref{fig:stability} shows that only the graded layer is accurate uniformly over the
sampled range $z\in[10^{-3},10^{12}]$, with $\mathcal E_S(1)=\univEone$, whereas the Legendre, quasi-uniform, and random
$\tanh$ layers give $\mathcal E_S(1)=\certELeg$, $\certESphere$, and $\certERandom$, respectively.

The size of $\mathcal E_S(1)$ of the graded layer is determined by the density of the slopes for moderate $z$ and by
the largest slope for large $z$. With $\nu=8$, $s_{\max}=10^7$ gives $\mathcal E_S(1)=\designA$, and $s_{\max}=10^{11}$
gives $\designB$. With $\nu=10$, we obtain $\designC$ for the truncation tolerance $10^{-12}$ of the singular values and
$\univEone$ for $10^{-16}$, and with $\nu=12$, $\designD$ with $r+1=\designDdim$. Hence, the density, the range, and the dimension of the compressed space have to be increased together, which requires the extended precision in
\cref{rem:implementation}.

\Cref{tab:inside} compares discretizations of one window with at most $\subUnknowns$ unknowns in time. At all three
times, the graded layer is more accurate than the polynomial, quasi-uniform, and random $\tanh$ layers, which have
errors of order $10^{-3}$ to $10^{-1}$ at $t-t_{n-1}=10^{-3}k$ inside the initial layer. The table also shows
that the geometric distribution of the slopes is essential, whereas the activation function is not. With slopes drawn
at random from the log-uniform distribution on the same range, the error at $10^{-3}k$ is $\ablRandomInside$ instead of
$\univInside$, with equally spaced slopes, the compressed space collapses to $28$ dimensions, and with the exponentials
$1-e^{-2s_i(\tau+1)}$ instead of $\tanh(s_i(\tau+1))$, the errors are even slightly smaller. Finally, compared with the
most accurate of $\hpSearched$ geometric refinements with at most $\subUnknowns$ unknowns (ratio, number of levels, and
degrees) for each time, the graded layer is more accurate at $10^{-3}k$ ($\univInside$ versus $\hpOptInside$ for the
ratio $\hpOptRatio$ with $\hpOptLevels$ levels and degrees from $\hpOptQa$ to $\hpOptQb$), whereas the geometric
refinement is far more accurate at $10^{-1}k$ and at the end of the window ($\hpOptMid$ and $\hpOptEnd$). The graded layer thus pays off deep inside the initial layer, at the price of the accuracy for the smooth part of the window.

\sbox{\tablebox}{\tablebodyfont
\begin{tabular}{@{}lrccc@{}}
\toprule
 & & \multicolumn{3}{c}{$(t-t_{n-1})/k$}\\
\cmidrule(lr){3-5}
one window & unknowns & $10^{-3}$ & $10^{-1}$ & $1$\\
\midrule
tanh, graded & 119 & \ensuremath{9.2{\times}10^{-10}} & \ensuremath{8.1{\times}10^{-11}} & \ensuremath{5.6{\times}10^{-11}} \\
Legendre $\mathbb{P}_{118}$ & 119 & \ensuremath{8.8{\times}10^{-3}} & \ensuremath{3.3{\times}10^{-4}} & \ensuremath{1.5{\times}10^{-3}} \\
tanh, quasi-uniform, $m{=}200$ & 83 & \ensuremath{1.3{\times}10^{-1}} & \ensuremath{2.1{\times}10^{-3}} & \ensuremath{2.0{\times}10^{-3}} \\
tanh, random, $m{=}192$ & 46 & \ensuremath{1.2{\times}10^{-1}} & \ensuremath{2.7{\times}10^{-3}} & \ensuremath{6.8{\times}10^{-3}} \\
tanh, random slopes & 112 & \ensuremath{4.1{\times}10^{-8}} & \ensuremath{2.0{\times}10^{-8}} & \ensuremath{3.0{\times}10^{-9}} \\
tanh, equally spaced slopes & 28 & \ensuremath{8.8{\times}10^{-1}} & \ensuremath{1.6{\times}10^{-1}} & \ensuremath{9.8{\times}10^{-2}} \\
$1-e^{-2s(\tau+1)}$, graded & 116 & \ensuremath{1.9{\times}10^{-10}} & \ensuremath{4.0{\times}10^{-11}} & \ensuremath{2.3{\times}10^{-11}} \\
$hp$, optimized for each time$^{*}$ & $\le119$ & \ensuremath{2.0{\times}10^{-8}} & \ensuremath{1.1{\times}10^{-13}} & \ensuremath{<10^{-13}} \\
\bottomrule
\end{tabular}}
\begin{table}[tbhp]
\caption{Worst-case modal errors $\mathcal E_S(1,\tau)=\sup_{z>0}|U_z(\tau)-y_z(\tau)|$ (\cref{thm:modal}) at three
times in one window of length $k$ with at most $\subUnknowns$ unknowns in time ($m$: number of neurons; the quasi-uniform
and random layers saturate at $83$ and $46$). $^*$For each time, the most accurate of $\hpSearched$ geometric
refinements.}\label{tab:inside}
\usetable
\end{table}

\subsection{Heat Equation with Rough Data}\label{sec:ex3}
\begin{figure}[tbhp]
\centering\includegraphics[width=\textwidth]{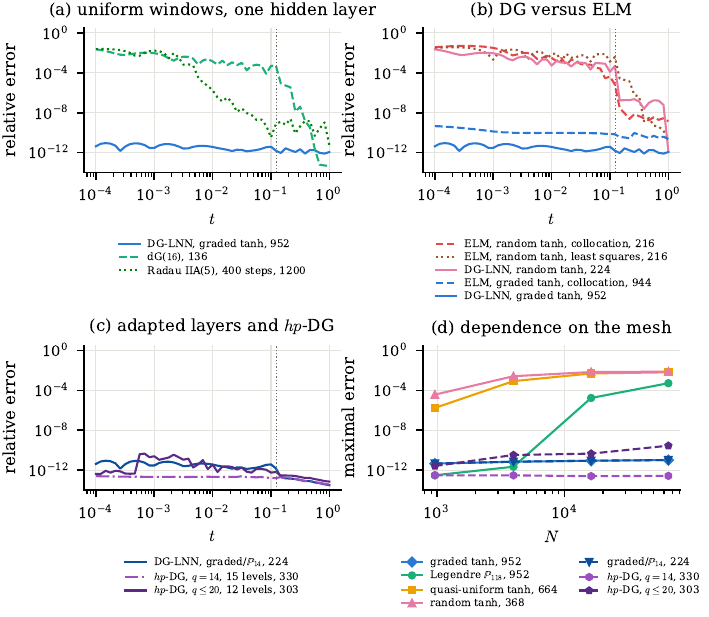}
\caption{Heat equation with rough data. (a)--(c) Relative errors at $49$ times in $[10^{-4},1]$ on a $128^2$ mesh, where
DG-LNN, the ELMs, and dG($16$) use $8$ uniform windows of length $k=1/8$, $hp$-DG refines the first of them, and the
Radau IIA(5) method uses $400$ uniform steps. (d) Maximal relative error over these times on meshes from $32^2$ to
$256^2$ (quasi-uniform: $m=200$, random: $m=192$). The numbers in the legends are the total numbers of unknowns in
time.}\label{fig:heat}
\end{figure}

We solve $u_t-0.05\Delta u=g(x)(1+\sin8\pi t)$ on the unit square with homogeneous Dirichlet boundary conditions and
$T=1$, where the initial value is the indicator function of the disk with center $(0.35,0.4)$ and radius $0.2$, and
$g(x)=10\exp(-|x-(0.65,0.6)|^2/0.01)$. \Cref{fig:heat}(a) shows the errors on a $128^2$ mesh ($N=\heatN$) at $49$ times
$t\in[10^{-4},1]$ equally spaced in $\log t$, where dG($16$) and DG-LNN use $8$ uniform windows of length $k=1/8$, and
the Radau IIA(5) method uses $400$ uniform steps. In the first window, which contains the initial layer, dG($16$) and the Radau IIA(5) method have errors up to $\timeLegFirst$ and
$\timeRadauFirst$, respectively, whereas DG-LNN with the graded layer is accurate uniformly in time, with the maximal
error $\timeGradedMax$, but it uses $\timeGradedUnk$ unknowns in time.

\Cref{fig:heat}(b) separates the effect of the hidden layer from that of the DG condition. With the random $\tanh$
layer, the ELMs with collocation and least squares as well as DG-LNN have errors of order $10^{-2}$ in the first
window, since this layer does not resolve the stiff modes. With the graded neurons, the ELM with
collocation resolves the initial layer, but its error $\timeELMGrMax$ is larger than that of DG-LNN by more than an
order of magnitude. Least squares with these neurons fails, since its matrix is numerically singular at the
Chebyshev--Gauss--Lobatto points (condition number $\lsCond$ at $z=10^4$), where the steep neurons vanish at all points
but $t_{n-1}$, and its amplification factor exceeds one ($\lsMaxR$) with points graded towards $t_{n-1}$. For decaying
modes, the hidden layer thus determines whether the initial layer is resolved, whereas the DG condition improves the
accuracy and makes every hidden layer stable.

\Cref{fig:heat}(c) compares the adapted layers with $hp$-DG. With graded/$\Pq_{\hpcQ}$, the maximal error is $\timeAdPMax$, the same as for the graded layer on all windows, with
$\timeAdPUnk$ unknowns in time, and graded/truncated gives $\timeAdMax$ with $\timeAdUnk$ unknowns. In agreement with
the worst-case modal errors in \cref{sec:setup}, $hp$-DG with $q=\hpcQ$ is more accurate ($\timeHPcMax$) but needs
$\timeHPcUnk$ unknowns, whereas the a priori plan gives $\timeAprMax$ with $\timeAprUnk$ unknowns.

\Cref{fig:heat}(d) shows the dependence on the mesh. As $N$ increases from $\meshNmin$ to $\meshNmax$, the largest
eigenvalue grows from $\meshLamFirst$ to $\meshLamLast$, and the error of the Legendre layer $\Pq_{\univr}$, which has
the size of the graded layer, increases from $\meshLegFirst$ to $\meshLegLast$, whereas the quasi-uniform and random
$\tanh$ layers reach $\meshSphereLast$ and $\meshRandomLast$, respectively. In contrast, the errors of the graded layer,
of graded/$\Pq_{\hpcQ}$, and of $hp$-DG with $q=\hpcQ$ remain at most $\meshGradedErr$, $\meshAdPErr$, and $\meshHPcErr$,
respectively, on all meshes, whereas the error of the a priori plan increases from $\meshAprFirst$ to $\meshAprLast$,
since its largest error occurs at $t\approx c/\lambda_{\max}$ in sub-windows with low degrees. Thus, a hidden layer or a geometric refinement that is chosen in advance for a range of $\lambda k$ by the worst-case modal error does not
deteriorate as the mesh is refined, in agreement with \cref{thm:modal}, which covers the part of the solution due to
the initial value.

In terms of computing time, \cref{tab:cost} lists the most accurate run of each method over the times
$t\in\{0.01,0.02,0.1,0.5,1\}$ on the $128^2$ mesh on one process. DG-LNN with the graded layer requires $r+1=\subUnknowns$
shifted factorizations, which take $\heatSetupDG$\,s and are reused in every window, and two sets of solves per window (the Schur solve and one step of refinement). With $W=8$, it reaches $\heatMaxDGeight$ in $\heatTimeDGeight$\,s,
whereas the Radau IIA(5) method with $3200$ steps reaches $\heatMaxbestRadau$ in $\heatMaxbesttimeRadau$\,s, and
dG($q$) on uniform windows does not reach $10^{-4}$. For moderate accuracies, however, it is not efficient, since its
factorizations alone take longer than the adaptive ESDIRK method with the tolerance $10^{-8}$, which reaches
$\heatESDIRKeight$ in $\heatESDIRKeighttime$\,s. The adapted layers graded/$\Pq_{\hpcQ}$ and graded/truncated reach
$\heatMaxAdPeight$ and $\heatMaxAdeight$ in $\heatTimeAdPeight$\,s and $\heatTimeAdeight$\,s with $\heatFactAdP$ and
$\heatFactAd$ factorizations, respectively, whereas $hp$-DG with $q=14$ reaches $\heatMaxHPcfourteen$ in
$\heatTimeHPcfourteen$\,s with $\heatFactHPcfourteen$ factorizations, and with $q=12$, $\heatMaxHPctwelve$ in
$\heatTimeHPctwelve$\,s with $\heatFactHPctwelve$ factorizations. $hp$-DG with the same degree as graded/$\Pq_{\hpcQ}$ is slower, since it requires one set of
factorizations for each length of the sub-windows, whereas with $q=12$, it requires fewer factorizations and is faster
and more accurate.

\sbox{\tablebox}{\tablebodyfont
\begin{tabular}{@{}llcrrrr@{}}
\toprule
method & run & error & factoriz. & solves & setup [s] & total [s]\\
\midrule
CN & steps 6400 & \ensuremath{4.4{\times}10^{-6}} & 1 & 6400 & -- & 4.86 \\
Radau IIA(5) & steps 3200 & \ensuremath{5.1{\times}10^{-11}} & 2 & 6400 & 0.12 & 15.49 \\
ELM, random $\tanh$ (collocation) & $W$ 32 & \ensuremath{4.1{\times}10^{-4}} & 27 & 864 & -- & 3.87 \\
DG-LNN, Legendre (dG) & $W$ 16 & \ensuremath{4.9{\times}10^{-4}} & 17 & 272 & 1.00 & 1.71 \\
DG-LNN, graded $\tanh$ & $W$ 8 & \ensuremath{3.8{\times}10^{-12}} & 119 & 1904 & 6.75 & 12.23 \\
DG-LNN, graded $\tanh$ & $W$ 16 & \ensuremath{2.0{\times}10^{-12}} & 119 & 3808 & 6.71 & 17.66 \\
\midrule
ESDIRK5 (adaptive) & tolerance $10^{-12}$ & \ensuremath{5.7{\times}10^{-11}} & 392 & 2744 & -- & 6.03 \\
$hp$-DG, linear degrees & $(q,L)$ (20,12) & \ensuremath{5.3{\times}10^{-12}} & 142 & 527 & -- & 9.70 \\
$hp$-DG, constant degree & $(q,L)$ (12,15) & \ensuremath{2.1{\times}10^{-12}} & 105 & 308 & -- & 6.91 \\
$hp$-DG, constant degree & $(q,L)$ (14,15) & \ensuremath{2.1{\times}10^{-13}} & 225 & 660 & -- & 14.90 \\
DG-LNN, graded/$\mathbb{P}_{14}$ & $W$ 8 & \ensuremath{3.8{\times}10^{-12}} & 134 & 448 & 7.55 & 8.87 \\
DG-LNN, graded/$\mathbb{P}_{14}$ & $W$ 16 & \ensuremath{2.0{\times}10^{-12}} & 134 & 688 & 7.54 & 9.50 \\
DG-LNN, graded/truncated & $W$ 8 & \ensuremath{3.8{\times}10^{-12}} & 159 & 798 & 8.95 & 11.26 \\
\bottomrule
\end{tabular}}
\begin{table}[tbhp]
\caption{Heat equation on the $128^2$ mesh on one process, most accurate run of each method (graded layer and
graded/$\Pq_{\hpcQ}$: $W=8$ and $16$; $hp$-DG: $q=12$ and $14$): maximal relative error over
$t\in\{0.01,0.02,0.1,0.5,1\}$, sparse factorizations and solves, time of the factorizations, and total time (median of
three runs). Top: one space on uniform steps or windows; bottom: adaptive steps, refined windows, or adapted layers.}\label{tab:cost}
\usetable
\end{table}

Since the shifted solves \eqref{eq:shifted} of a window are independent, they can be distributed over $P$ processes, for
which we use the diagonalized solver with $\ell_{\rm ref}$ refinement steps set by $\kappa(\mathcal S)$. The same holds
for $hp$-DG, whose sub-windows, however, are solved in sequence. Hence, each process first factorizes its share of the
shifts of all windows and sub-windows, and then the windows are solved in sequence, each by $1+\ell_{\rm ref}$ rounds of
independent shifted solves. \Cref{tab:parallel} shows the wall-clock times on a $192^2$ mesh ($N=\parN$) with $8$
windows (medians of three runs, which differ by at most $\parSpread\%$). The graded layer on all windows requires
$\ell_{\rm ref}=\parRounds$ because of $\kappa(\mathcal S)=\parCondS$, which gives $\parGStages$ rounds, and it takes
$\parGTmax$\,s on $\parPmax$ processes with the error $\parGErr$, whereas the Schur solver with one refinement step takes $\parSchurT$\,s on one process with the
error $\parSchurErr$. The adapted layers graded/$\Pq_{\hpcQ}$ require $\parAdPStages$ rounds, compared with
$\parHPStages$ for $hp$-DG, whose $\parHPWindows$ windows and sub-windows require $\parHPRoundsPer$ rounds each. They
take $\parAdPTone$\,s on one process and $\parAdPTmax$\,s on $\parPmax$ processes, compared with $\parHPTone$\,s and
$\parHPTmax$\,s for $hp$-DG, whereas $hp$-DG is more accurate ($\parHPErr$ compared with $\parAdPErr$). The difference is
in the factorizations, which take $\parAdPSetup$\,s and $\parHPSetup$\,s on $\parPmax$ processes, respectively, whereas
the time stepping after them takes $\parAdPLoop$\,s and $\parHPLoop$\,s, so that the fewer rounds of graded/$\Pq_{\hpcQ}$
do not make its time stepping faster on these processes. In fact, the first window of the graded layer requires one set
of $\parAdFirstShifts$ shifted factorizations, of which $\parAdFirstReal$ are real, whereas that of $hp$-DG requires one
set of $\parHPShifts$, of which $\parHPReal$ is real, for each of its $\parHPSubLengths$ lengths of the sub-windows, and
its window loop is not accelerated beyond $\parHPShifts$ processes.

\sbox{\tablebox}{\tablebodyfont
\begin{tabular}{@{}lrrrrrrrrc@{}}
\toprule
method & unknowns & factoriz. & rounds & $P=1$ & $P=2$ & $P=4$ & $P=8$ & $P=16$ & error\\
\midrule
DG-LNN, graded $\tanh$ & 952 & 91 & 56 & 35.39 & 20.70 & 13.05 & 9.64 & 7.66 & \ensuremath{8.1{\times}10^{-11}} \\
DG-LNN, graded/truncated & 399 & 111 & 28 & 20.33 & 11.13 & 6.45 & 4.23 & 3.10 & \ensuremath{5.5{\times}10^{-12}} \\
DG-LNN, graded/$\mathbb{P}_{14}$ & 224 & 99 & 28 & 15.65 & 8.60 & 4.90 & 3.31 & 2.48 & \ensuremath{5.5{\times}10^{-12}} \\
$hp$-DG, $q=14$, $L=15$ & 330 & 120 & 66 & 23.82 & 13.13 & 7.15 & 4.11 & 3.05 & \ensuremath{2.5{\times}10^{-13}} \\
\bottomrule
\end{tabular}}
\begin{table}[tbhp]
\caption{Parallel runs on a $192^2$ mesh with $8$ windows and the diagonalized solver: unknowns in time,
factorizations, sequential rounds of independent shifted solves, wall-clock time in seconds on $P$ processes (median of
three runs, including all factorizations), and maximal relative error over $t\in\{0.01,0.02,0.1,0.5,1\}$.}\label{tab:parallel}
\usetable
\end{table}

\subsection{Heat Equation on Unstructured Meshes}\label{sec:plate}
\begin{figure}[tbhp]
\centering\includegraphics[width=\textwidth]{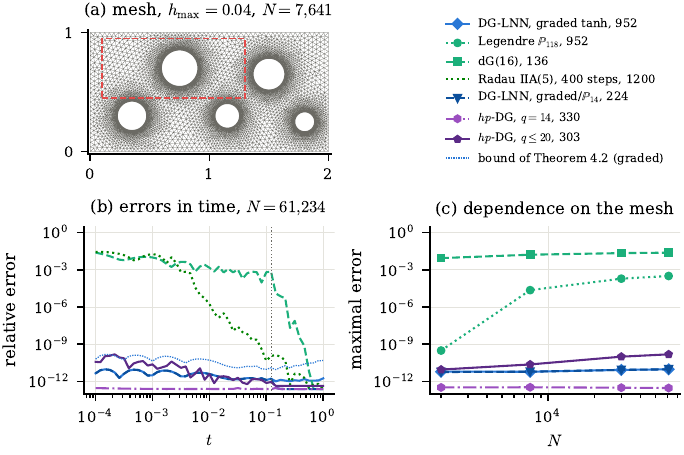}
\caption{Heat equation on a perforated plate. (a) Mesh with $h_{\max}=0.04$ and support of the initial value
(dashed). (b) Relative errors at $49$ times in $[10^{-4},1]$ on the finest mesh and the bound of \cref{thm:modal} for
the graded layer (dotted). (c) Maximal relative error over these times on the four meshes. The numbers in the legend,
which is shared by (b) and (c), are the total numbers of unknowns in time.}\label{fig:plate}
\end{figure}

Next, we consider unstructured, locally refined meshes with mixed boundary conditions, which the finite element space
imposes exactly. This test concerns the time stepping for the stiffer spectra of locally refined meshes, not neural
network discretizations in space. We solve $u_t-0.05\Delta u=0$ for $t\in(0,1]$ on the
plate $(0,2)\times(0,1)$ with five circular holes, where $u=0$ on the boundaries of the holes, $\partial_nu=0$ on the
outer boundary, and the initial value is the indicator function of $[0.1,1.3]\times[0.45,0.95]$ (\cref{fig:plate}(a)).
The meshes are generated by Gmsh \cite{GeuzaineRemacle2009} with the element size $h_{\max}$ away
from the holes and $h_{\max}/8$ on the holes for $h_{\max}=0.08$, $0.04$, $0.02$, and $0.014$, which gives $N$ from $\plateNmin$ to $\plateNmax$, and the largest eigenvalue reaches $\plateLamLast$ on the finest mesh, compared with
$\meshLamLast$ on the finest uniform mesh in \cref{sec:ex3}. Since $f=0$ and $a$ is symmetric and coercive,
\cref{thm:modal} applies. On the meshes with at most $8000$ unknowns, the reference differs from the exact
solution by at most $\plateRefExact$.

\Cref{fig:plate}(b) shows the errors on the finest mesh with uniform windows of length $k=1/8$ and $400$ steps of the
Radau IIA(5) method. As on the uniform meshes, dG($16$) and the Radau IIA(5) method have errors up to
$\plateFirstWindG$ and $\plateFirstWinRadau$, respectively, in the first window, whereas the error of DG-LNN with the
graded layer is at most $\plateGradedFirst$. The bound of \cref{thm:modal}, $\mathcal E_S(n,\tau)$ with the supremum over the eigenvalues of this mesh, $z\in[\lambda_{\min}k,\lambda_{\max}k]$, relative to $\norm{u_h(t)}$, lies
above the error at all sample times, by factors from $\plateBoundMin$ to $\plateBoundMax$. The bound is attained for the
eigenfunctions: on the coarsest mesh, for the initial value $u_0=w_\ell$, where $w_\ell$ is the eigenfunction for which
the modal term at $t=10^{-3}k$ is the largest ($\lambda_\ell k\approx\sharpZ$), the error at this time differs from the
modal term $\sharpTerm$ by the relative amount $\sharpDiff$. \Cref{fig:plate}(c) shows that the maximal errors of the
graded layer, of graded/$\Pq_{\hpcQ}$, and of $hp$-DG with $q=\hpcQ$ remain at most $\plateGradedErr$, $\plateAdPErr$,
and $\plateHPcErr$, respectively, on all meshes, whereas those of the Legendre layer $\Pq_{\univr}$ with the same number
of unknowns and of the a priori plan increase from $\plateLegFirst$ to $\plateLegLast$ and from $\plateAprFirst$ to
$\plateAprLast$, respectively. Finally, we refine the meshes on the holes by the factor $32$ instead of $8$
($h_{\max}=0.02$, $N=\stiffN$) and use one window of length $k=1$, so that $\lambda_{\max}k=\stiffLk$. Then the maximal
error of the graded layer is $\stiffGraded$, whereas $\Pq_{\univr}$ and dG($16$) give $\stiffLeg$ and $\stiffdG$,
respectively. Thus, the graded layer, which is fixed in advance, remains accurate when local refinement increases
$\lambda_{\max}k$.

\subsection{High-Order Elements}\label{sec:highorder}
\begin{figure}[tbhp]
\centering\includegraphics[width=\textwidth]{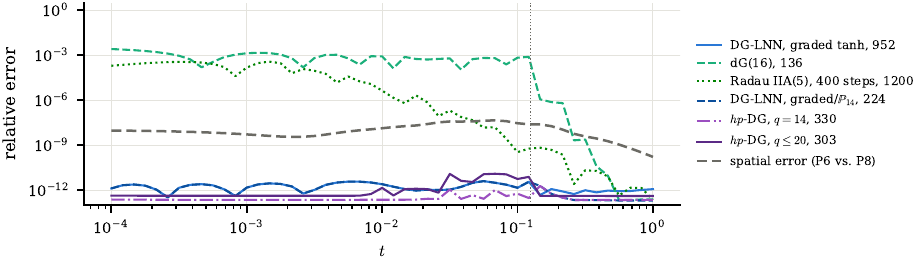}
\caption{Heat equation with P\hoDeg{} elements: relative errors of the time discretizations at $49$ times in
$[10^{-4},1]$ (windows of length $k=1/8$, the first one refined by $hp$-DG) and the spatial error, estimated by P8
elements. The numbers in the legend are the total numbers of unknowns in time.}\label{fig:highorder}
\end{figure}

The time discretization errors in \cref{sec:ex3,sec:plate} are far below the spatial errors of P1 elements. Hence, we
also consider a regime where the spatial error is small. We solve $u_t-0.05\Delta u=0$ on the unit square with $u=0$
on the boundary and $T=1$, where the initial value $u_0=\exp(-|x-(0.45,0.55)|^2/0.03^2)$ is a narrow Gaussian. We use
P\hoDeg{} elements on an unstructured mesh with the element size $0.012$ near the center of $u_0$ and $0.1$ away from it ($N=\hoN$), so that $\lambda_{\max}k=\hoLk$ for $k=1/8$. The spatial error, estimated by the
difference to P8 elements on the same mesh ($N=\hoNhi$), is between $\hoSpaceMin$ and $\hoSpaceMax$ at the
sample times (\cref{fig:highorder}). In the first window, dG($16$) and the Radau IIA(5) method with $400$ steps have errors up to $\hodGFirst$ and $\hoRadauFirst$, respectively, which exceed the spatial error by orders of
magnitude. In contrast, the maximal error of the graded layer is $\hoGradedMax$, and those of
graded/$\Pq_{\hpcQ}$, $hp$-DG with $q=\hpcQ$, and the a priori plan are $\hoAdPMax$, $\hoHPcMax$, and $\hoAprMax$,
respectively. With high-order elements, a time stepping method therefore has to resolve the initial layer in order not to
spoil the spatial accuracy, and the graded layer does so with the same hidden layer on uniform windows.

\subsection{Long-Time Wave Propagation}\label{sec:ex4}
\begin{figure}[tbhp]
\centering\includegraphics[width=\textwidth]{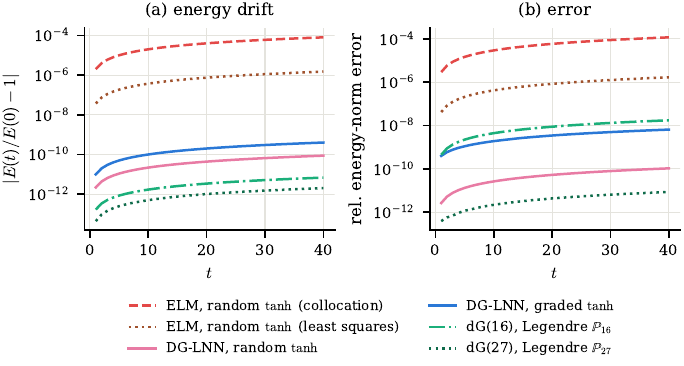}
\caption{Energy drift and energy-norm error for the wave equation at the ends of every tenth window. The ELMs and
DG-LNN use the same random $\tanh$ layer with $32$ neurons.}\label{fig:wave}
\end{figure}

Finally, we consider the wave equation, which is conservative and has no initial layer. We solve $u_{tt}=\Delta u$ on
the unit square with $u=0$ on the boundary, $u(0)=\exp(-|x-(0.5,0.45)|^2/0.004)$, and $u_t(0)=0$, discretized by P1
elements on a $48^2$ mesh, in the first-order form
$\diag(K,M)\,x'+\left(\begin{smallmatrix}0&-K\\K&0\end{smallmatrix}\right)x=0$ with $x=(u,u_t)$, whose operator is skew-symmetric in the energy inner product, which replaces $\ip\cdot\cdot$ in this section. Since $a$
is not coercive, \eqref{eq:dg} is uniquely solvable by \cref{thm:spectral}(ii) because $\sigma(D_n)\cap i\R=\emptyset$
for the hidden layers used here. We use $400$ windows of length $0.1$ up to $T=40$ with $\waveN$ spatial unknowns.
\Cref{fig:wave} shows the energy drift and the error at the ends of every tenth window. Per window, the ELMs have $27$
unknowns in time, DG-LNN with the random layer $28$, with the graded layer $\subUnknowns$, and dG($16$) and dG($27$)
have $17$ and $28$. The ELMs with collocation and with least squares have the energy drifts \waveelmE{} and
\waveelmlsE{} and the errors \waveelmerr{} and \waveelmlserr{} at $t=40$, where the amplification factor of least
squares exceeds one on the imaginary axis by up to \waveelmlsgrowth{} per window. In contrast, DG-LNN with the same
hidden layer gives \wavedgrandE{} and \wavedgranderr{} at a similar cost (\waveelmtime\,s, \waveelmlstime\,s, and \wavedgrandtime\,s, respectively).

These results do not depend on the particular random layer. Since DG-LNN, collocation, and least squares in the norm
induced by $M^{-1}$ all decouple in the eigenvectors of $(K,M)$ into the scalar problems $y'=i\omega y$, we compute their
errors mode by mode for $\wmSeeds$ random layers. On the $48^2$ mesh ($\omega_{\max}k\approx\wmOmegaA$), the errors at
$t=40$ have the medians $\wmDGAmed$, $\wmCollAmed$, and $\wmLSAmed$ for DG-LNN, collocation, and least squares, and the
largest errors over the layers are $\wmDGAhi$, $\wmCollAhi$, and $\wmLSAhi$, respectively. On the $\wmMeshB^2$ mesh
($\omega_{\max}k\approx\wmOmegaB$), collocation diverges for $\wmCollBdiv$ of the $\wmSeeds$ layers, and for
$\wmRadBdiv$ with the Radau points, since its amplification factor has poles near the imaginary axis, and least squares
has the median $\wmLSBmed$, whereas the errors of DG-LNN remain at most $\wmDGBhi$. The DG condition thus improves the median accuracy of the same hidden layer by three (least squares) to five (collocation) orders of magnitude, and it
is stable for every layer. DG-LNN with the graded layer is also stable, and its error \wavedggradederr{} lies between
those of dG($16$) (\wavedglegerr) and dG($27$) (\wavedglegTwentySevenerr), so that it is less efficient than polynomials for this problem, since its neurons are designed for decaying modes.
\FloatBarrier
\section{Conclusion}\label{sec:conclusion}
In this paper, we used the hidden layer of a shallow neural network as the trial space of DG time stepping for finite
element discretizations. DG-LNN keeps the A- and L-stability of DG time stepping for every hidden layer, and on uniform
windows, its error for symmetric, coercive, homogeneous problems is bounded for all spatial discretizations by the
worst-case modal error of the hidden layer. Since this quantity is computed once on a reference window, it can be used
to choose hidden layers in advance, as we did for the graded layer, the adapted layers, and $hp$-DG. The analysis also explains why polynomials, which attain the largest order
for non-stiff modes, do not resolve stiff initial layers on uniform windows, since the reproducing kernel determines how
strongly unresolved modes are damped.

The experiments separate the roles of the hidden layer and of the DG condition. For decaying modes, the hidden layer decides whether a stiff
initial layer is resolved. With geometrically graded slopes, DG-LNN reached about
$10^{-11}$ on uniform windows on all meshes tested, whereas polynomial, quasi-uniform, and random $\tanh$ layers of the
same size did not, and the distribution of the slopes mattered more than the activation function. The DG condition, in
turn, makes a given hidden layer reliable. With the graded neurons, collocation was less accurate and least squares
failed, and in the long-time wave test, the DG condition improved the median accuracy of random hidden layers by three
to five orders of magnitude.

The graded layer is not a uniform improvement over polynomials. Used on the first window instead of the sub-windows of
$hp$-DG with $q=\hpcQ$, it saved factorizations and was faster on one and on $\parPmax$ processes, but $hp$-DG was more
accurate, its time stepping after the factorizations was faster, and with $q=12$, it was faster on one process. For
moderate accuracy, an adaptive Runge--Kutta method was cheaper, and for oscillating modes, polynomials were more
efficient. We therefore expect the graded layer to be most useful when high accuracy is required inside stiff initial
layers and the factorizations are a large part of the cost.

The limitations of the present work are as follows. First, the analysis and the experiments are restricted to linear
problems with time-independent coefficients, and the experiments to two space dimensions with at most $\stiffN$ spatial
unknowns, sparse direct solvers, and $\parPmax$ processes of one machine. Second, the matrices in time of the graded
layer have to be formed in extended precision, and they are large and far from normal ($\norm A\approx\gradedNormA$,
$\kappa(\mathcal S)=\refCondS$), so that the solves in time require iterative refinement, and the diagonalized solver is
less accurate than the Schur form. Third, except for the test with high-order elements, the time discretization errors
considered here are far below the spatial discretization errors. Finally, the error bounds cover neither forcing terms
nor the wave equation. As future work, we will extend DG-LNN to nonlinear problems using window-wise Newton--Krylov
solvers and to large-scale problems in three space dimensions, where the reuse of the factorizations is more important,
and we will study a bound of the worst-case modal error of the graded layer in terms of its slopes.

\appendix
\section{Computational Details}\label{sec:certificates}
In this section, we describe how $\mathcal E_S(n)$ and $\mathcal E_S(n,\tau)$ in \cref{sec:ex1} are computed. We use
$R_S(-z)=\phi(1)^\top(A+\tfrac z2G)^{-1}\phi(-1)$ by \cref{thm:spectral}(iii) and
$U_z(\tau)=\phi(\tau)^\top(A+\tfrac z2G)^{-1}\phi(-1)$ with the matrices in time of \cref{rem:implementation}, whose symmetric part is exact, and for $m$ sub-windows, the factors of the sub-windows are multiplied. The suprema over $z$
are taken on a logarithmic grid of $700$ points in $[10^{-3},10^{12}]$, and values below $10^{-13}$ are at the rounding
level of evaluating $R_S^n$. For A-stability, the closed left half-plane is sampled at the points $iy$ with $500$
logarithmically spaced $y\in[10^{-3},10^4]$ and at the points $\rho e^{i\theta}$ with $80$ logarithmically spaced
$\rho\in[10^{-2},10^7]$ and $41$ equispaced $\theta\in[\pi/2,3\pi/2]$.

For the graded layer and for the layers in \cref{tab:inside} with other slopes or exponentials, the matrices are formed
with $\grDps$ digits by the Gauss--Legendre rule with $32$ points on each interval between the breakpoints $-1+j/4$,
$j=0,\dots,8$, and $-1+10^{-13}\cdot1.3^j<1$, $j=0,1,\dots$, and the check in \cref{rem:implementation} uses
$\grCheckDps$ digits, $40$ points, and the breakpoints $-1+10^{-16}\cdot1.25^j$. For the other hidden layers, we use
double precision and the Gauss--Legendre rule with $20$ points on each interval between the
breakpoints $-1+j/4$, $j=0,\dots,8$, and $-1+\delta$, where $\delta$ takes $60$ values equally spaced in $\log\delta$
between $10^{-10}$ and $2$, and the compression in \cref{rem:implementation} is the singular value decomposition,
without weights, of the matrix of the functions $\sigma_i(w_i\tau+b_i)-\sigma_i(-w_i+b_i)$ at $\max(64,16m)$ Chebyshev
points and at the points $-1+\delta$ for as many values of $\delta$ equally spaced in $\log\delta$ between $10^{-9}$ and
$2$, where $m$ is the number of neurons.

\bibliographystyle{siamplain}
\bibliography{references}

\end{document}